\documentclass[12pt]{article}
\usepackage{geometry}
\usepackage{titlesec}

\titleformat{\section}
  {\normalfont\fontsize{14}{15}\bfseries}{\thesection}{1em}{}
\titleformat{\subsection}
{\normalfont\fontsize{14}{15}\bfseries}{\thesubsection}{1em}{}

\newcommand{\footremember}[2]{%
   \footnote{#2}
    \newcounter{#1}
    \setcounter{#1}{\value{footnote}}%
}
\providecommand{\keywords}[1]
{
  \textbf{\textit{Keywords---}} #1
}

\usepackage{graphicx}
\usepackage{caption}
\usepackage{subcaption}%
\usepackage{amsmath,amsthm}      
\usepackage{color}
\usepackage{latexsym}
\usepackage{hyperref}
\numberwithin{equation}{section}
\newtheorem{defn}{Definition}
\newtheorem{theorem}{Theorem}

\begin{document}

\title{Entropy stable flux correction for scalar hyperbolic conservation laws 
}


\author{Sergii Kivva\footremember{alley}{Institute of Mathematical Machines and System Problems, National Academy of Sciences, Ukraine}
}



\date{\Large}

\maketitle

\begin{abstract}
It is known that Flux Corrected Transport algorithms can produce entropy-violating solutions of hyperbolic conservation laws.
Our purpose is to design flux correction with {\it maximal} antidiffusive fluxes to obtain entropy solutions of scalar hyperbolic conservation laws. To do this we consider a hybrid difference scheme that is a linear combination of a monotone scheme and a scheme of high-order accuracy. Flux limiters for the hybrid scheme are calculated from a corresponding optimization problem. Constraints for the optimization problem consist of inequalities that are valid for the monotone scheme and applied to the hybrid scheme. 
We apply the discrete cell entropy inequality with the proper numerical entropy flux to single out a physically relevant solution of scalar hyperbolic conservation laws. A nontrivial approximate solution of the optimization problem yields expressions to compute the required flux limiters.
We present examples that show that not all numerical entropy fluxes guarantee to single out a physically correct solution of scalar hyperbolic conservation laws.

\keywords{entropy solution, flux corrected transport, scalar conservation laws, difference scheme, linear programming}
\end{abstract}

\section{Introduction}
\label{intro}
We consider the numerical solution of the Cauchy problem for 1D scalar hyperbolic conservation laws
\begin{eqnarray} 
\label{eq:1} 
&& \frac{\partial u}{\partial t}+\frac{\partial }{\partial x}f(u) = 0  
 \qquad t > 0 
 \qquad {-\infty  < x < \infty}  \\
 &&  u(x,0) = u_{0}(x)
\label{eq:2} 
\end{eqnarray}
where $f$ is a smooth flux function.

It is well known \cite{b1} that solutions of \eqref{eq:1}-\eqref{eq:2} may develop singularities at a finite time even for a smooth initial condition. Hence, we should interpret \eqref{eq:1} in the sense of distribution and search for weak solutions. However, such weak solutions are not unique. To single out a unique physically relevant weak solution, the latter should satisfy
\begin{equation}
\label{eq:3} 
\frac{\partial U(u)}{\partial t} + \frac{\partial F(u)}{\partial x} \le 0	
\end{equation}
in the sense of distribution for every entropy pair ($U,F$). Here $U$ is a convex function, the so-called entropy function and its entropy flux $F$ such that satisfy $F'(u) = U'(u)f'(u)$.

We discretize \eqref{eq:1} by the conservative difference scheme 
\begin{equation}
\label{eq:4} 
\frac{1}{\Delta t}({\hat y_i} - {y_i}) + \frac{1}{\Delta x}\left[ {{h_{i + {1 / 2}}} - {h_{i - {1 / 2}}}} \right] = 0 								
\end{equation}
where the numerical flux ${h_{i + {1 / 2}}}$ is calculated by
  
\begin{equation}
\label{eq:5} 
{h_{i + {1 / 2}}} = h_{i + {1 / 2}}^L + {\alpha _{i + {1 / 2}}}\left[ {h_{i + {1 / 2}}^H - h_{i + {1 / 2}}^L} \right] 								
\end{equation}
Here, ${y_i} = y({x_i},t)$ is the approximation value at the gridpoint (${x_i} = i\Delta x,t)$; ${\hat y_i} = y({x_i},t + \Delta t)$; $\Delta x $ and $\Delta t $ are the spatial and temporal computational grid size, respectively. $h_{i + {1 / 2}}^H$ and $h_{i + {1 / 2}}^L$ are a high-order and low-order numerical fluxes such that ${h_{i + {1 / 2}}} = h({y_{i - l + 1}},...,{y_{i + r}})$ is the Lipschitz continuous numerical flux consistent with the differential flux, that is $h(u,...,u) = f(u)$ for all flux-limiters ${\alpha _{i + {1 / 2}}} \in \left[ {0,1} \right]$.
 
The expression in square brackets on the right-hand side of \eqref{eq:5} can be considered as an antidiffusive flux. For flux-correction, we will calculate the flux-limiters ${\alpha _{i + {1 / 2}}}$ by using linear programming. Flux corrected transport (FCT) was firstly developed by Oran and Boris~\cite{b2} for solving the transient continuity equation. Later Zalesak~\cite{b4,b3} extended FCT to multidimensional explicit difference schemes. Several implicit FEM-FCT schemes for unstructured grids were proposed by Kuzmin and his coworkers~cite{b5,b6,b8,b7}. However, known FCT algorithms do not guarantee to obtain entropy solutions for hyperbolic conservation laws.
 
We discretize the entropy inequality \eqref{eq:3} as follows 
\begin{equation} 
 U({\hat y_i}) - U({y_i}) + \frac{{\Delta t}}{{\Delta x}}\left[ {{H_{i + {1  / 2}}} - {H_{i - {1 / 2}}}} \right] \le 0 
 \label{eq:6} 
\end{equation}
where ${H_{i + {1 / 2}}} = H({y_{i - l + 1}},...,{y_{i + r}})$ is the numerical entropy flux consistent with the differential one $H(u,...,u) = F(u)$. 

A difference scheme \eqref{eq:4} is called {\it entropy stable} if computed solutions satisfy the discrete cell entropy inequality \eqref{eq:6}. We mention here the pioneering studies of entropy stable schemes by Lax~\cite{b37}. 
Entropy stable schemes were developed by several authors \cite{b43,b16,b18,b17,b39}. Harten et al.~\cite{b43} showed that for scalar conservation laws all explicit monotone schemes are entropy stable. A class of $E$-schemes which includes monotone schemes and is entropy stable was introduced by Osher~\cite{b16,b18}. 
In \cite{b12,b13}~Tadmor studied the entropy stability of difference approximations to nonlinear hyperbolic conservation laws by comparing the entropy production of a given scheme against properly chosen entropy-conservative schemes.
Chalons and LeFloch~\cite{b40} studied limiting solutions of fully discrete finite-difference schemes for a diffusive-dispersive conservation law. They investigated the dependence of these scheme solutions on balance between dissipative and disperse forces and producing non-classical shock solutions violating the standard entropy criterion. In \cite{b42}~LeFloch et al. proposed a general approach to construct second and third-order accurate, fully discrete implicit entropy-conservative schemes. Construction of other high-order accurate, fully discrete entropy stable schemes can be found in \cite{b36,b38,b41,b35}.

To single out a physically relevant solution, we use the so-called proper numerical entropy flux concept of which was formulated by Merriam~\cite{b9} and Sonar~\cite{b10}. Zhao and Wu~\cite{b11} proved that a three-point monotone semi-discrete schemes in conservative form satisfy corresponding semi-discrete entropy inequality with proper numerical entropy flux. The numerical entropy flux $H({y_{i - l + 1}},...,{y_{i + r}})$ for $F$ is not unique. Tadmor~\cite{b36,b12} proposed another form of the numerical entropy flux and introduced general families of entropy-conservative schemes. 
Note that not all numerical entropy fluxes guarantee that discrete cell entropy inequality \eqref{eq:6} only holds for a physically relevant solution. Numerical examples show that the numerical entropy flux proposed by Tadmor \cite{b36,b12} is one of them. 
 
In this paper, we construct the flux correction for 1D scalar hyperbolic conservation laws to obtain numerical entropy solutions for which antidiffusive fluxes are maximal. For this, flux limiters for the hybrid scheme \eqref{eq:4}-\eqref{eq:5} are computed from the optimization problem with constraints that are valid to the low-order scheme. Approximate solution of the optimization problem yields the desired flux correction formulas.

The paper is organized as follows: In Section \ref{sec:1}, we present the estimations that are valid for a monotone scheme and the definition of the proper numerical entropy flux. The optimization problem for finding flux limiters and the algorithm of its solving are described in Section \ref{sec:2}. The approximate solution of the optimization problem is derived in Section \ref{sec:3}. The results of numerical experiments are given in Section \ref{sec:4}. Concluding remarks are drawn in Section \ref{sec:5}.

\section{Monotone Difference Scheme}
\label{sec:1}

We consider a three point low-order scheme \eqref{eq:4} in the form
\begin{equation}
\label{eq:21}
{\hat y_i} - {y_i} + \frac{{\Delta t}}{{\Delta x}}\left[ {h_{i + {1 /2}}^L({y_i},{y_{i + 1}}) - h_{i - {1 /2}}^L({y_{i - 1}},{y_i})} \right] = 0	\end{equation}
where the low-order numerical flux $h_{i+{1/2}}^L = h(v,u)$ is consistent and satisfies the following inequalities
\begin{equation} 
\label{eq:22}
\frac{\partial h(v,u)}{\partial v} \ge 0  \qquad
\frac{\partial h(v,u)}{\partial u} \le 0  \qquad
\forall \quad u,v\in R
\end{equation}

\begin{theorem}
\label{th1}
 Let the  numerical fluxes $h_{i + {1 /2}}^L$ of three point scheme \eqref{eq:21} satisfy the inequalities \eqref{eq:22} for all ${y_i}$. If $\Delta t$ satisfies
\begin{equation}
\label{eq:23}
\Delta t \mathop {\max }\limits_i \left[ {\frac{{\partial h_{i+1/2}^L}}{\partial {y_i}} - \frac{{\partial h_{i-1/2}^L}}{\partial {y_i}}} \right] \le \Delta x 										
\end{equation}
then the scheme \eqref{eq:21} is monotone and the following inequality holds
\begin{equation}
\label{eq:24}
\begin{split} 
\min ({y_{i - 1}},{y_i},{y_{i + 1}}) & \le {y_i} - \frac{{\Delta t}}{{\Delta x}}\left[ {h_{i + {1 /2}}^L({y_i},{y_{i + 1}}) - h_{i - {1 /2}}^L({y_{i - 1}},{y_i})} \right]  \\
& \le \max ({y_{i - 1}},{y_i},{y_{i + 1}}) 	\end{split}			  						
\end{equation}
\end{theorem}

\begin{proof}
 The monotonicity of the scheme \eqref{eq:21} follows from the conditions of the theorem.

Consider the function 
\[ \Phi ({y_{i - 1}},{y_i},{y_{i + 1}}) = {y_i} - \frac{\Delta t}{\Delta x} \left[ {h_{i + {1 /2}}^L({y_i},{y_{i + 1}}) - h_{i - {1 /2}}^L({y_{i - 1}},{y_i})} \right] \] 
that is non-decreasing for its arguments. Denote by ${\bar y_i} = \max ({y_{i - 1}},{y_i},{y_{i + 1}})$ and ${\underline y_i} = \min ({y_{i - 1}},{y_i},{y_{i + 1}})$. Then we have
\begin{equation*}
  {\underline y}_i = \Phi ({\underline y}_i,{\underline  y_i},{\underline y_i }) \le 
\Phi ({y_{i - 1}},{y_i},{y_{i + 1}})  \le \Phi ({\bar y}_i,{\bar  y_i},{\bar y_i })  =  {\bar y_i} 
\end{equation*}
This completes the proof of the theorem.
\end{proof}

For the conservative scheme \eqref{eq:4} we define the proper numerical entropy flux in a similar way to Zhao and Wu \cite{b11}.

\begin{defn}
 The numerical entropy flux $H({y_{i - l + 1}},...,{y_{i + r}})$ of the conservative scheme \eqref{eq:4} is called proper if
\begin{equation}
\label{eq:26}  
\begin{split} 
& \frac{\partial }{{\partial {y_j}}}H({y_{i - l + 1}},...,{y_{i + r}}) = \\  & \qquad =\frac{{dU({y_j})}}{{d{y_j}}}\frac{\partial }{{\partial {y_j}}}h({y_{i - l + 1}},...,{y_{i + r}}) \quad for \quad j = i - l + 1,...,i + r  
\end{split} 		 
\end{equation} 
\end{defn}

Then the proper numerical entropy flux for the difference scheme \eqref{eq:4} and \eqref{eq:5} can be written in the form

\begin{equation}
\label{eq:27}  
  {H_{i + {1 / 2}}} = H_{i + {1 / 2}}^L + {\alpha _{i + {1 / 2}}}\left( {H_{i + {1 / 2}}^H - H_{i + {1 / 2}}^L} \right) 							\end{equation}
where $H_{i+{1/2}}^L$ and $H_{i+{1/2}}^H$ are the low-order and high-order proper numerical entropy fluxes corresponding to the numerical fluxes $h_{i+{1/2}}^L$ and $h_{i+{1/2}}^H$.

\begin{theorem}
\label{th2}
Assume that $f\in{C^1}(R)$, a convex function $U\in{C^2}(R)$ and for all $i$ $\Delta t$ satisfies the inequalities  
\begin{equation}
\label{eq:28a}  
 \Delta t\left[ {\frac{{\partial h_{i+ {1/2}}^L}}{{\partial {y_i}}} - \frac{{\partial h_{i - {1 /2}}^L}}{{\partial {y_i}}}} \right] \le \Delta x	
 	\end{equation} 
 
\begin{equation}
\label{eq:28b}  
 \Delta t\mathop {\max }\limits_{s \in [{{\underline y}_i},{\bar y_i}]} U''(s)\;{\left( {h_{i + {1 /2}}^L - h_{i - {1 /2}}^L} \right)^2} \le 2\Delta x\left[ {U'({y_i})\left( {h_{i + {1 /2}}^L - h_{i - {1 /2}}^L} \right) - H_{i + {1 /2}}^L + H_{i - {1 /2}}^L} \right] 		  		
\end{equation}
where ${{\underline y}_i} = \min ({y_{i - 1}},{y_i},{y_{i + 1}})$  and ${\bar y_i} = \max ({y_{i - 1}},{y_i},{y_{i + 1}})$.

Then the fully discrete monotone scheme \eqref{eq:21} satisfies the discrete cell entropy inequality
\begin{equation}
\label{eq:28c}  
U({\hat y_i}) - U({y_i}) + \frac{{\Delta t}}{{\Delta x}}\left[ {H_{i + {1/2}}^L({y_i},{y_{i + 1}}) - H_{i - {1/2}}^L({y_{i - 1}},{y_i})} \right] \le 0		\end{equation}
where $H_{i + {1/2}}^L$ is the proper numerical entropy flux corresponding to the numerical flux $h_{i + {1/2}}^L$.
\end{theorem}

\begin{proof}
Multiplying the equation \eqref{eq:21} by $U'({y_i})$ and subtracting it from the left-hand side of the inequality \eqref{eq:28c}, we have
\begin{equation}
\label{eq:28d}  
\begin{split}  
  U({\hat y_i}) - U({y_i}) \; + & \; \frac{{\Delta t}}{{\Delta x}}\left[ {H_{i + {1 /2}}^L - H_{i - {1 /2}}^L} \right] \\
 =  U({\hat y_i}) - U({y_i}) -  U'({y_i})\left( {{{\hat y}_i} - {y_i}} \right)  \; + & \; \frac{{\Delta t}}{{\Delta x}}\left[ {H_{i + {1 /2}}^L - H_{i - {1 /2}}^L - U'({y_i})\left( {h_{i + {1/2}}^L - h_{i - {1 /2}}^L} \right)} \right]  \\
 = \frac{1}{2}{\left( {\frac{{\Delta t}}{{\Delta x}}} \right)^2}  U''(\xi ){({\hat y_i} - {y_i})^2} \; + & \; \frac{{\Delta t}}{{\Delta x}}\left[ {\int\limits_{{y_i}}^{{y_{i + 1}}} {(U'(v) - U'({y_i}))\frac{\partial }{{\partial v}}h({y_i},v)dv} }\right. \\
  + & \; \left. { \int\limits_{{y_{i - 1}}}^{{y_i}} {(U'(v) - U'({y_i}))\frac{\partial }{{\partial v}}h(v,{y_i})dv} } \right]
\end{split}
\end{equation}
where $\xi  = \theta {\hat y_i} + (1 - \theta ){y_i}$, $0 < \theta  < 1$. 
 
It is easy to see that the first term on the right-hand side of \eqref{eq:28d} is non-negative, and the second and third terms in square brackets are non-positive. Indeed, for a convex function $U(w) \in {C^1}(R)$, the inequality $\left[ {U'(w) - U'(v)} \right](w - v) \ge 0$ holds for any $w,v \in R$. Therefore, the sign of the expression $\left[ {U'(w) - U'(v)} \right]$ coincides with ${\mathop{\rm sgn}} (w - v)$ and the integrands do not change the sign on the integration interval. Hence,
\[{\mathop{\rm sgn}} \left[ {\int\limits_{{y_i}}^{{y_{i + 1}}} {\left( {U'(v) - U'({y_i})} \right)} \frac{\partial }{{\partial v}}h({y_i},v)dv} \right] =  - {{\mathop{\rm sgn}} ^2}({y_{i + 1}} - {y_i})\;\]						
\[{\mathop{\rm sgn}} \left[ {\int\limits_{{y_{i - 1}}}^{{y_i}} {\left( {U'(v) - U'({y_i})} \right)} \frac{\partial }{{\partial v}}h(v,{y_i})dv} \right] =  - {{\mathop{\rm sgn}} ^2}({y_i} - {y_{i - 1}})\;\]						
This completes the proof of the theorem.

\end{proof}

Thus, to obtain a physically relevant solution of \eqref{eq:1} and \eqref{eq:2} the flux limiters $\alpha _{i+{1/2}}$ should satisfy 

\begin{equation}
\label{eq:28}  
  {\hat y_i} - {y_i} + \frac{{\Delta t}}{{\Delta x}}\left[ {h_{i + {1 / 2}}^L + {\alpha _{i + {1 / 2}}}h_{i + {1 / 2}}^{AD} - h_{i - {1 / 2}}^L - {\alpha _{i - {1 / 2}}}h_{i - {1 / 2}}^{AD}} \right] = 0  					\end{equation}
\begin{equation}
\label{eq:29}  
  U({\hat y_i}) - U({y_i}) + \frac{{\Delta t}}{{\Delta x}}\left[ {H_{i + {1 / 2}}^L + {\alpha _{i + {1 / 2}}}H_{i + {1 / 2}}^{AD} - H_{i - {1 / 2}}^L - {\alpha _{i - {1 / 2}}}H_{i - {1 / 2}}^{AD}} \right] \le 0   		\end{equation}
where $h_{i + {1/2}}^{AD} = h_{i + {1/2}}^H - h_{i + {1/2}}^L$ and $H_{i + {1/2}}^{AD} = H_{i + {1/2}}^H - H_{i + {1/2}}^L$.

\section{Finding Flux Limiters}
\label{sec:2}
For the numerical solution of the Cauchy problem \eqref{eq:1}-\eqref{eq:2}, the modeling area is replaced by the large but limited area. Therefore, in further, we suggest that the flux limiters are finding in the limited area 
\[{U_{ad}} = \left\{ {\left. \boldsymbol{\alpha}  \right| \quad {\rm{   0}} \le {\alpha _{i + {1 / 2}}} \le {\rm{1,  \quad   - }}N \le i \le N} \right\}\]
To find the flux limiters ${\alpha _{i + {1 / 2}}}$, we consider the following optimization problem
\begin{equation}
\label{eq:31} 
\Im (\boldsymbol{\alpha} ) = \sum\limits_i {{c_{i + {1 / 2}}}} {\alpha _{i + {1 / 2}}} \to \mathop {\max }\limits_{\boldsymbol{\alpha}  \in {U_{ad}}}  									
\end{equation}
subject to 
\begin{multline}
\label{eq:32}  
 \frac{{\Delta x}}{{\Delta t}}\left( {\mathop {\min }\limits_{j \in {S_i}} {y_j} - {y_i}} 
 \right)  + h_{i + {1 / 2}}^{L} - h_{i - {1/2}}^{L} \le  - {\alpha _{i + {1 / 2}}}\;
  h_{i + {1/2}}^{AD} + {\alpha _{i-{1/2}}}\; h_{i - {1 / 2}}^{AD} \le 	\\		
  \le \frac{{\Delta x}}{{\Delta t}}\left( {\mathop {\max }\limits_{j \in {S_i}} {y_j} - {y_i}} \right) + h_{i + {1/2}}^{L} - h_{i - {1/2}}^{L} 
\end{multline} 

\begin{multline}
\label{eq:33}  
 \frac{{\Delta x}}{{\Delta t}}\left[ {U({{\hat y}_i}) - U({y_i}) - {b_i}\left( {{{\hat y}_i} -  {y_i}} \right)} \right] + H_{i + {1 / 2}}^{L} - {b_i}h_{i + {1 / 2}}^{L} - H_{i - {1 / 2}}^{L} + {b_i}h_{i - {1 / 2}}^{L} \le \\				
  \le {\alpha _{i + {1 / 2}}}\;({b_i} h_{i + {1 / 2}}^{AD} - H_{i + {1 / 2}}^{AD}) - {\alpha _{i - {1 / 2}}}\;({b_i} h_{i - {1 / 2}}^{AD} - H_{i - {1 / 2}}^{AD})
\end{multline}

\begin{equation}
\label{eq:34}
 \frac{{\Delta x}}{{\Delta t}}\left( {{{\hat y}_i} - {y_i}} \right) + h_{i + {1 / 2}}^{L} - h_{i - {1 / 2}}^{L} + {\alpha _{i + {1 / 2}}} h_{i + {1 / 2}}^{AD} - {\alpha _{i - {1 / 2}}} h_{i - {1 / 2}}^{AD} = 0					
\end{equation}
where  ${\textit{\textbf{ c}}} \ge 0$ and $\textbf{\textit{ b}}$ are a-priori specified vectors, $S_i$ is set of indexes $\left\{ i-1,i,i+1 \right\} $.

Due to constraints \eqref{eq:33}, the optimization problem \eqref{eq:31}-\eqref{eq:34} is nonlinear. Therefore, finding the numerical entropy solution of the Cauchy problem \eqref{eq:1}-\eqref{eq:2} at one time step can be represented as the following iterative process.

\begin{description}
\item[\it Step 1.] Specify positive numbers $\delta $, $\epsilon_1 $, $\epsilon_2 $ and  vectors $ \textit{\textbf{c}}$, $\textbf{\textit{b}}$. Set $p=0$, 
 $ {\boldsymbol{\hat y}}^0 = \boldsymbol{y}$, $ \boldsymbol{\alpha} ^0 = 0$.
\item[\it Step 2.] Find $ \boldsymbol {\alpha}^{p + 1}$ as the solution of the following linear programming problem
\begin{equation}
\label{eq:35}
\Im (\boldsymbol {\alpha} ) = \sum\limits_i {{c_{i + {1 / 2}}}} \alpha _{i + {1 / 2}}^{p + 1} \to \mathop {\max }\limits_{\boldsymbol {\alpha} ^{p + 1} \in {U_{ad}}} 						
\end{equation}
subject to
\begin{multline}
\label{eq:36}
 \frac{{\Delta x}}{{\Delta t}}\left( {\mathop {\min }\limits_{j \in {S_i}} {y_j} - {y_i}} \right)  + h_{i + {1 / 2}}^{L} - h_{i - {1 / 2}}^{L} \le  - \alpha _{i + {1/2}}^{p+1}\; h_{i + {1/2}}^{AD} + \alpha _{i - {1 / 2}}^{p + 1}\; h_{i - {1 / 2}}^{AD} \le 		\\		
  \le \frac{{\Delta x}}{{\Delta t}}\left( {\mathop {\max }\limits_{j \in {S_i}} {y_j} - {y_i}} \right) + h_{i + {1/2}}^{L} - h_{i - {1/2}}^{L} 
\end{multline}

\begin{multline}
\label{eq:37}
\frac{{\Delta x}}{{\Delta t}}\left[ {U(\hat y_i^p) - U({y_i}) - {b_i}\left( {\hat y_i^p - {y_i}} \right)} \right] + H_{i + {1/2}}^{L} - {b_i}h_{i + {1/2}}^{L} - H_{i - {1/2}}^{L} + {b_i}h_{i - {1/2}}^{L} \le \\			
 \le \alpha _{i + {1/2}}^{p+1} \left( {b_i} h_{i + {1/2}}^{AD} - H_{i + {1/2}}^{AD} \right) - \alpha _{i - {1/2}}^{p + 1} \left( {b_i} h_{i - {1/2}}^{AD} - H_{i - {1/2}}^{AD} \right)
\end{multline}

\item[\it Step 3.] For the ${\boldsymbol {\alpha} ^{p + 1}}$, find ${\boldsymbol {\hat y}^{p + 1}}$ from the system of linear equations
\begin{equation}
\label{eq:38}
 \frac{{\Delta x}}{{\Delta t}}\left( {\hat y_i^{p + 1} - {y_i}} \right) + h_{i + {1/2}}^{L} - h_{i - {1/2}}^{L} + \alpha _{i + {1/2}}^{p + 1} h_{i + {1 / 2}}^{AD} - \alpha _{i - {1/2}}^{p + 1} h_{i - {1/2}}^{AD} = 0					
\end{equation}

\item[\it Step 4.] Algorithm stop criterion
\begin{equation}
\label{eq:39}
 \frac{{\left| {\hat y_i^{p + 1} - \hat y_i^p} \right|}}{{\max \left( {\delta ,\left| {\hat y_i^{p + 1}} \right|} \right)}} < {\varepsilon _1},  \quad    \left| {\alpha _{i + {1 /2}}^{p + 1} - \alpha _{i + {1 / 2}}^p} \right| < {\varepsilon _2} 						
\end{equation}
If conditions \eqref{eq:39} hold, then $\boldsymbol {\hat y} = {\boldsymbol {\hat y}^{p + 1}}$. Otherwise, $p = p + 1$ and go to \textit{Step 2}.
\end{description}

\begin{theorem} 
\label{th3} 
 Let $f \in C^1(R)$ and $\Delta t$ satisfies the inequalities \eqref{eq:28a} and \eqref{eq:28b}, and $ \boldsymbol{\alpha}  \in {U_{ad}}$. 
Then the linear programming problem \eqref{eq:35}-\eqref{eq:37} is solvable.
\end{theorem}

\begin{proof}
 To prove the solvability of the linear programming problem \eqref{eq:35}-\eqref{eq:37}, it is sufficient to show that the objective function $\Im (\boldsymbol {\alpha} )$ is bounded on $U_{ad}$ and the feasible set is non-empty. The boundedness of the objective function \eqref{eq:35} follows from the boundedness of the vector $\boldsymbol {\alpha}$ on $U_{ad}$. For $\boldsymbol {\alpha}=0$ and under restriction \eqref{eq:28a}, the difference scheme \eqref{eq:28} is a three-point monotone for which the constraints \eqref{eq:36} hold. From Theorem \ref{th2} follows that under conditions \eqref{eq:28a} and \eqref{eq:28b} three-point monotone scheme satisfies the cell entropy inequalities \eqref{eq:37}. Therefore, the feasible set is nonempty and the linear programming problem \eqref{eq:35}-\eqref{eq:37} is solvable.
\end{proof}

\section{Approximate Solution of the Optimization Problem }
\label{sec:3}

We find a nontrivial $\boldsymbol {\alpha}  \in {U_{ad}}$ satisfying the inequalities \eqref{eq:36} which are rewritten in the form
\begin{equation}
\label{eq:41}
 - {\alpha _{i + {1 / 2}}} h_{i + {1/2}}^{AD} + {\alpha _{i - {1/2}}} h_{i - {1/2}}^{AD} \ge Q_i^ - 								
\end{equation}

\begin{equation}
\label{eq:42}
 - {\alpha _{i + {1/2}}} h_{i + {1/2}}^{AD} + {\alpha _{i - {1 / 2}}} h_{i - {1/2}}^{AD} \le Q_i^ + 								
\end{equation}
where 
\begin{eqnarray*}
  && Q_i^ +  = \frac{{\Delta x}}{{\Delta t}}\left( {\mathop {\max }\limits_{j \in {S_i}} {y_j} - {y_i}} \right) + \left[ {h_{i + {1 / 2}}^{L} - h_{i - {1 / 2}}^{L}} \right] \\
  && Q_i^ -  = \frac{{\Delta x}}{{\Delta t}}\left( {\mathop {\min }\limits_{j \in {S_i}} {y_j} - {y_i}} \right) + \left[ {h_{i + {1 / 2}}^{L} - h_{i - {1 / 2}}^{L}} \right]
\end{eqnarray*} 
Denote by $\alpha _i^ - $ and $\alpha _i^ + $ maximal values of components $\alpha$ for negative and positive terms on the left-hand side of \eqref{eq:41}-\eqref{eq:42}, respectively. Then
\begin{eqnarray}
\label{eq:43}
 - {\alpha _{i + {1 / 2}}} h_{i + {1/2}}^{AD} + {\alpha _{i - {1/2}}} h_{i - {1/2}}^{AD} \ge \alpha _i^ - P_i^ - \\							
 - {\alpha _{i + {1/2}}} h_{i + {1/2}}^{AD} + {\alpha _{i - {1/2}}} h_{i - {1/2}}^{AD} \le \alpha _i^ + P_i^ + 							
\label{eq:44}
\end{eqnarray}
where
\begin{eqnarray*}
  && P_i^ -  = \min \left( {0, -  h_{i + {1/2}}^{AD}} \right) + \min \left( {0, h_{i - {1/2}}^{AD}} \right)		\\				
  && P_i^ +  = \max \left( {0, -  h_{i + {1/2}}^{AD}} \right) + \max \left( {0, h_{i - {1/2}}^{AD}} \right) 
\end{eqnarray*} 
						
Each flux limiter ${\alpha _{i + {1 / 2}}}$ appears twice in \eqref{eq:43} and twice in \eqref{eq:44} with coefficients that differ only sign. Therefore ${\alpha _{i + {1 / 2}}}$ should not exceed 
\begin{equation}
\label{eq:45}
{\bar \alpha _{i + {1 / 2}}} = \left\{ {\begin{array}{*{20}{c}}
{\min (\alpha _i^ + ,\alpha _{i + 1}^ - ) = \min (R_i^ + ,R_{i + 1}^ - ), \quad {\rm{    }} h_{i + {1 / 2}}^{AD} < 0}\\
{\min (\alpha _i^ - ,\alpha _{i + 1}^ + ) = \min (R_i^ - ,R_{i + 1}^ + ), \quad {\rm{    }} h_{i + {1 / 2}}^{AD} > 0}
\end{array}} \right. 					
\end{equation} 
where $R_i^ \pm  = \min (1,{{Q_i^ \pm } \mathord{\left/
 {\vphantom {{Q_i^ \pm } {P_i^ \pm )}}} \right.
 \kern-\nulldelimiterspace} {P_i^ \pm )}}$. 
 
Note that the relations \eqref{eq:45} are similar to the FCT formulas of Zalesak \cite{b4} and Kuzmin \cite{b6}.

We rewrite \eqref{eq:37} in the form
\begin{equation}
\label{eq:46}
 {W_i} \le {\alpha _{i + {1/2}}}{d_{ii+1}} + {\alpha _{i - {1/2}}}{d_{ii - 1}}			\end{equation} 
where
\begin{eqnarray*}
 && {W_i} = \frac{{\Delta x}}{{\Delta t}}\left( {{{\hat U}_i} - {U_i} - {b_i}({{\hat y}_i} - {y_i})} \right) + H_{i + {1 / 2}}^{L} - {b_i}h_{i + {1 / 2}}^{L} - H_{i - {1 / 2}}^{L} + {b_i}h_{i - {1 / 2}}^{L} \\
 && {d_{ik}} = \left( {b_i} {h_{(i + k)/2}^{AD}} - {H_{(i + k)/2}^{AD}} \right) {\rm sgn} (k-i)
\end{eqnarray*} 
 
Reasoning similarly to above, from \eqref{eq:46} follows that upper bound of ${\alpha _{i + {1 / 2}}}$ is
\begin{equation}
\label{eq:47} 
\begin{split} 
  {\tilde \alpha _{i + {1/2}}} =  \min \left\{ 1, \frac{-W_i}{Y_i}\min \left( {0,{\mathop{\rm sgn}} \, d_{ii+1}} \right) +\max \left( {0, {\rm sgn} \, d_{ii+1}} \right), \right. \\
 \left. \frac{-W_{i + 1}}{Y_{i+1}} \min \left( {0,{\mathop{\rm sgn}} \, d_{i+1i}} \right) +\max \left( {0, {\rm sgn} \, d_{i+1i}} \right) 
  \right\}	 
\end{split}   				
\end{equation} 
where
${Y_i} = \min (0,{d_{ii+1}}) + \min (0,{d_{ii-1}})$ 

Thus, the nontrivial feasible solution of the linear programming problem \eqref{eq:35}-\eqref{eq:37} on ${U_{ad}}$ is equal to
\begin{equation}
\label{eq:48}
{\alpha _{i + {1 / 2}}} = \min ({\bar \alpha _{i + {1 / 2}}},{\tilde \alpha _{i + {1 / 2}}}) \end{equation}

\section{Numerical Examples}
\label{sec:4}

In this Section for the numerical solution of nonlinear scalar conservation laws \eqref{eq:1}, we apply a monotone scheme with the Rusanov numerical flux \cite{b14}

\begin{equation}
\label{eq:51} 
 h_{i + {1 / 2}}^{Rus}({y_i},{y_{i + 1}}) = \frac{1}{2}\left[ {f({y_i}) + f({y_{i + 1}}) - \mathop {\max }\limits_{s \in [{y_i},{y_{i + 1}}]} \left| {f'(s)} \right|({y_{i + 1}} - {y_i})} \right] 				 
\end{equation}  
If $f(u)$ is a convex function then $\mathop {\max }\limits_{s \in \left[ {{y_i},{y_{i + 1}}} \right]} \left| {f'(s)} \right|$ reduces to \[\mathop {\max }\limits_{s \in \left[ {{y_i},{y_{i + 1}}} \right]} \left| {f'(s)} \right|\; = \max \left( {\left| {f'({y_i})} \right|,\left| {f'({y_{i + 1}})} \right|} \right)\].

The proper numerical entropy flux for the Rusanov numerical flux \eqref{eq:51} can be written as
\begin{equation}
\label{eq:52} 
\begin{split} 
 H_{i + {1 / 2}}^{Rus}({y_i},{y_{i + 1}}) & = 
 \frac{1}{2}  \left[ \vphantom{\max \limits_{s \in [{y_i},{y_{i + 1}}]}} { F({y_i}) + F({y_{i + 1}}) -{}} \right. \\  
& -{} \left. \mathop {\max }\limits_{s \in [{y_i},{y_{i + 1}}]} \left| {f'(s)} \right| {\left(\mathstrut {U({y_{i + 1}}) - U({y_i})} \right)} \right]
\end{split} 		 
\end{equation} 

 In addition to the Rusanov flux \eqref{eq:51}, we use a monotone scheme with the Godunov numerical flux
\begin{equation}
\label{eq:53}
 f_{i + {1 / 2}}^G({y_i},{y_{i + 1}}) = \left\{ {\begin{array}{*{20}{c}}
{\mathop {\min \,}\limits_{{y_i} \le y \le {y_{i + 1}}} f(y) \quad {\rm{     if }} \; {y_i} \le {y_{i + 1}}}\\
{\mathop {\max \,}\limits_{{y_{i + 1}} \le y \le {y_i}} f(y) \quad {\rm{     if }} \; {y_i} > {y_{i + 1}}}
\end{array}} \right.	
\end{equation} 

Rusanov scheme and Godunov scheme are $E$-schemes \cite{b16,b17}. Osher \cite{b16,b18} showed that $E$ schemes satisfy the entropy inequality and have no less numerical viscosity than that of the Godunov scheme. Therefore in further the Godunov scheme is considered as reference scheme.

For solving linear programming problems we apply GLPK package v.4.65 $\href{url} {(https://www.gnu.org/software/glpk/)}$. 

\subsection{One-Dimensional Scalar Nonconvex Conservation Law}
\label{sec:41}

We consider the solution of the following Riemann problem \cite{b15}
\begin{eqnarray}
\label{eq:54}
 && \frac{\partial u}{\partial t} + \frac{\partial }{\partial x}f(u) = 0 		\\									
 && u(x,0) = 
  \begin{cases} 
   2 &    \text{if  } x \leq 1 \\
   -2 &   \text{if  } x > 1
  \end{cases}
\label{eq:55}
\end{eqnarray}
where 
\begin{equation}
\label{eq:56}
 f(u) = \frac{1}{4}({u^2} - 1)({u^2} - 4) 									
\end{equation} 

For the numerical solution of the IVP \eqref{eq:54}-\eqref{eq:56}, in the cell entropy inequality \eqref{eq:6} along with the proper numerical entropy flux \eqref{eq:26}, we will also use the numerical entropy flux proposed by Tadmor \cite{b12}
\begin{equation}
\label{eq:57}
 {H_{i + {1 / 2}}^T} = \frac{1}{2}({v_i} + {v_{i + 1}}) \, {h_{i + {1/2}}}(u(v)) - \frac{1}{2}[\psi ({v_i}) + \psi ({v_{i + 1}})] 					
\end{equation} 
where $v$ is the entropy variable $v = U'(u)$. If we assume that $U$ is strictly convex, then the mapping $u \to v$ is one-to-one and we can make the change of variables. $\psi $ is so-called entropy flux potential that is defined as
\begin{equation*}
  \psi (v) = vf(u(v)) - F(u(v))
\end{equation*} 
It is easy to see that the Tadmor's numerical entropy flux \eqref{eq:55} is not a proper entropy flux.

\begin{figure}[!t]
  \centering 
  \includegraphics[scale=0.95]{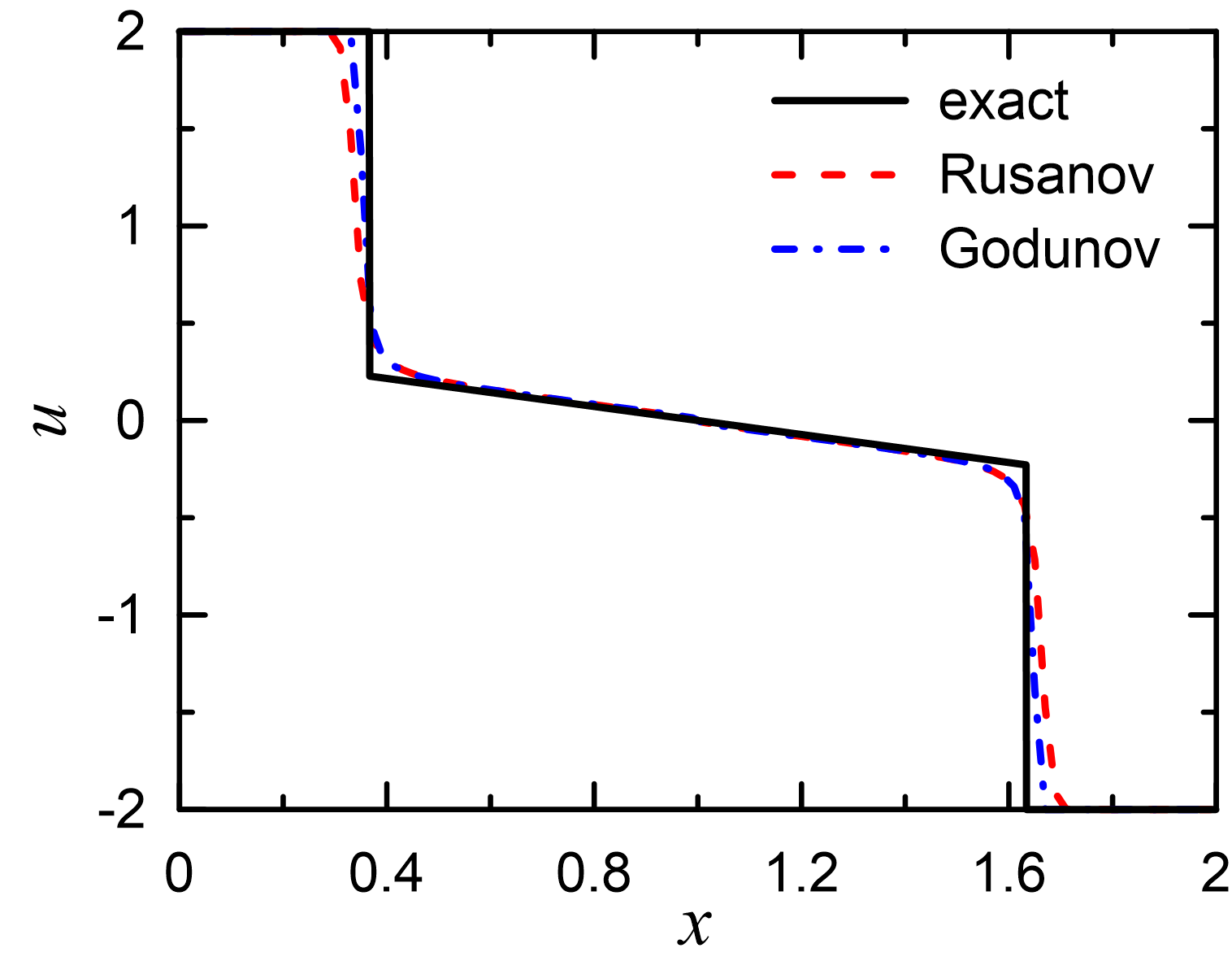}
\caption{Comparison between numerical solutions of the Riemann problem for the scalar nonconvex conservation law  \eqref{eq:54}-\eqref{eq:56} with the Godunov and Rusanov schemes}
\label{fig:1}       
\end{figure}
%
\begin{figure*}[!t]
  \includegraphics[width=0.5\textwidth]{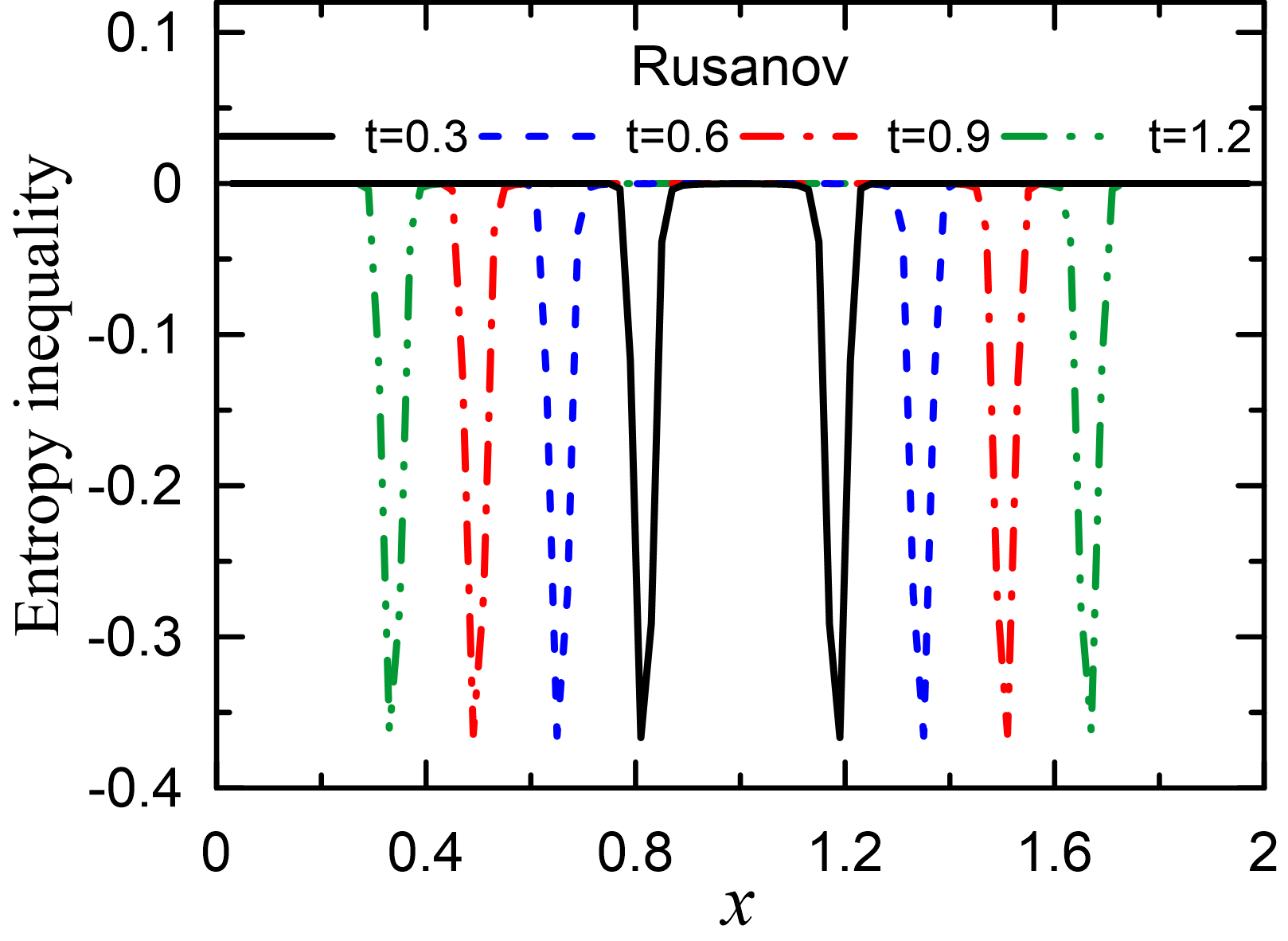}
  \includegraphics[width=0.5\textwidth]{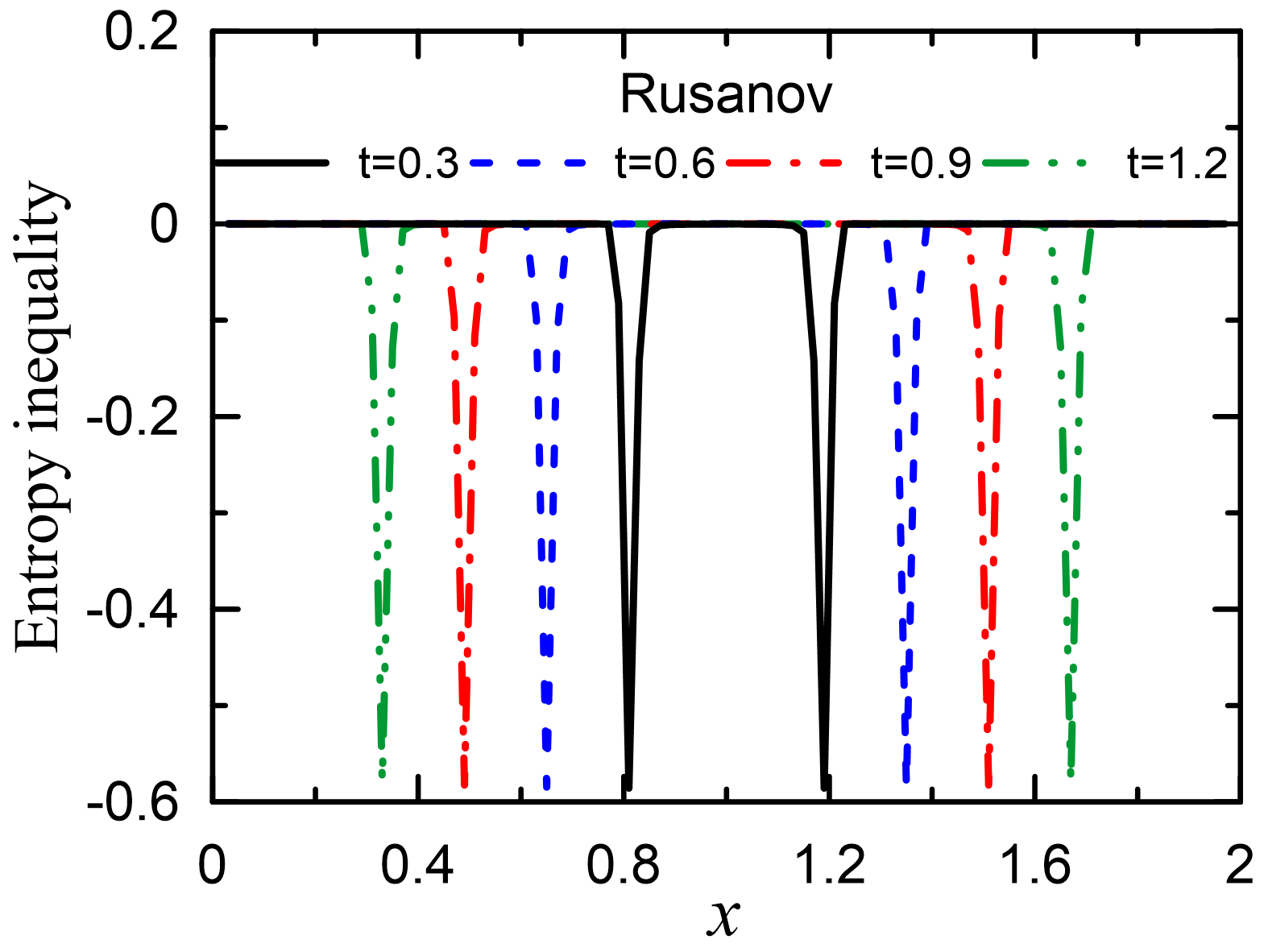}
\caption{Values of the discrete cell entropy inequality with the Tadmor's (left) and the proper (right) numerical entropy fluxes for the numerical solution of the Riemann problem for the scalar nonconvex conservation law  \eqref{eq:54}-\eqref{eq:56} with the Rusanov scheme at different times}
\label{fig:2}       
\end{figure*}
%

\begin{figure*}[!b]
  \includegraphics[width=0.5\textwidth]{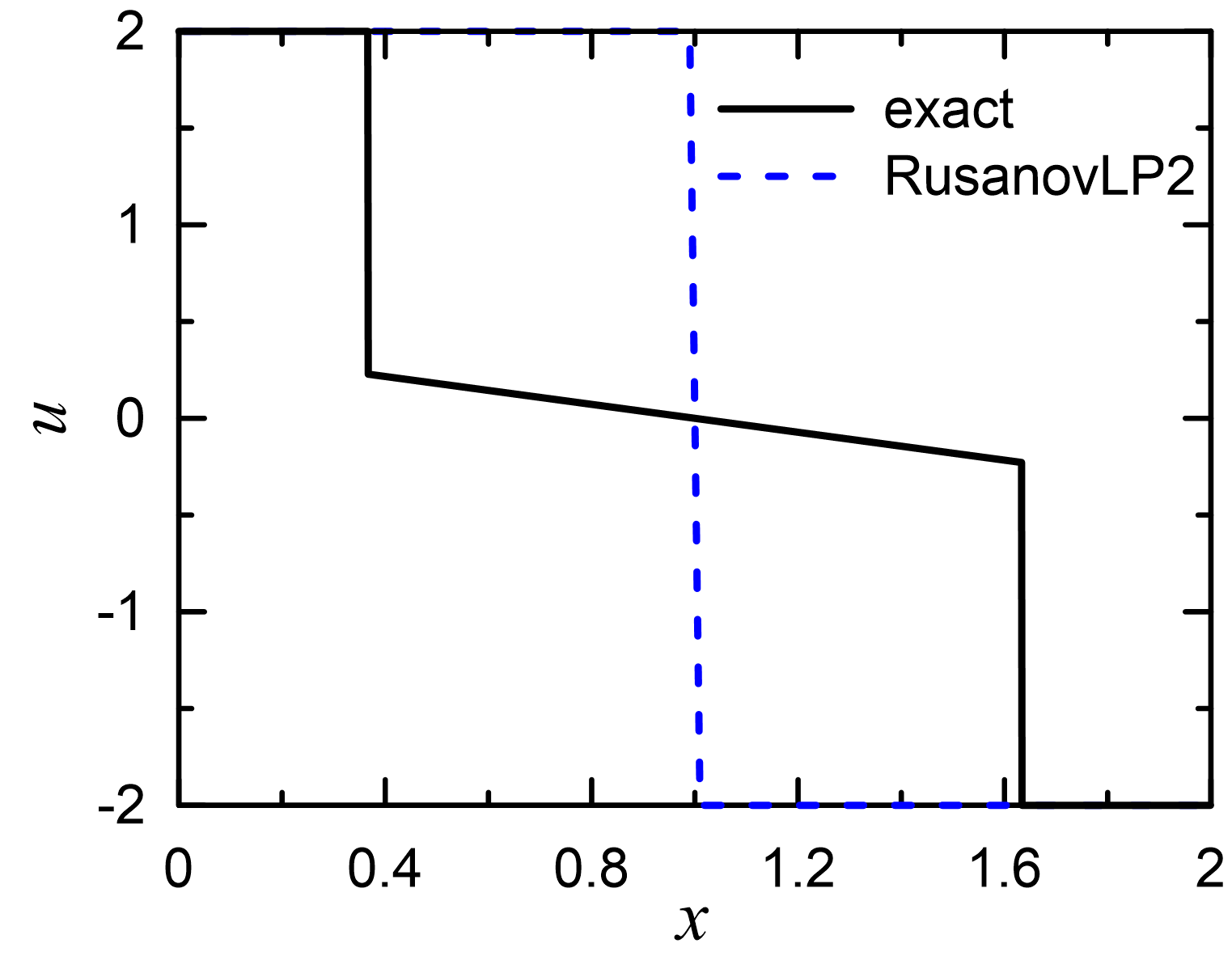}
  \includegraphics[width=0.5\textwidth]{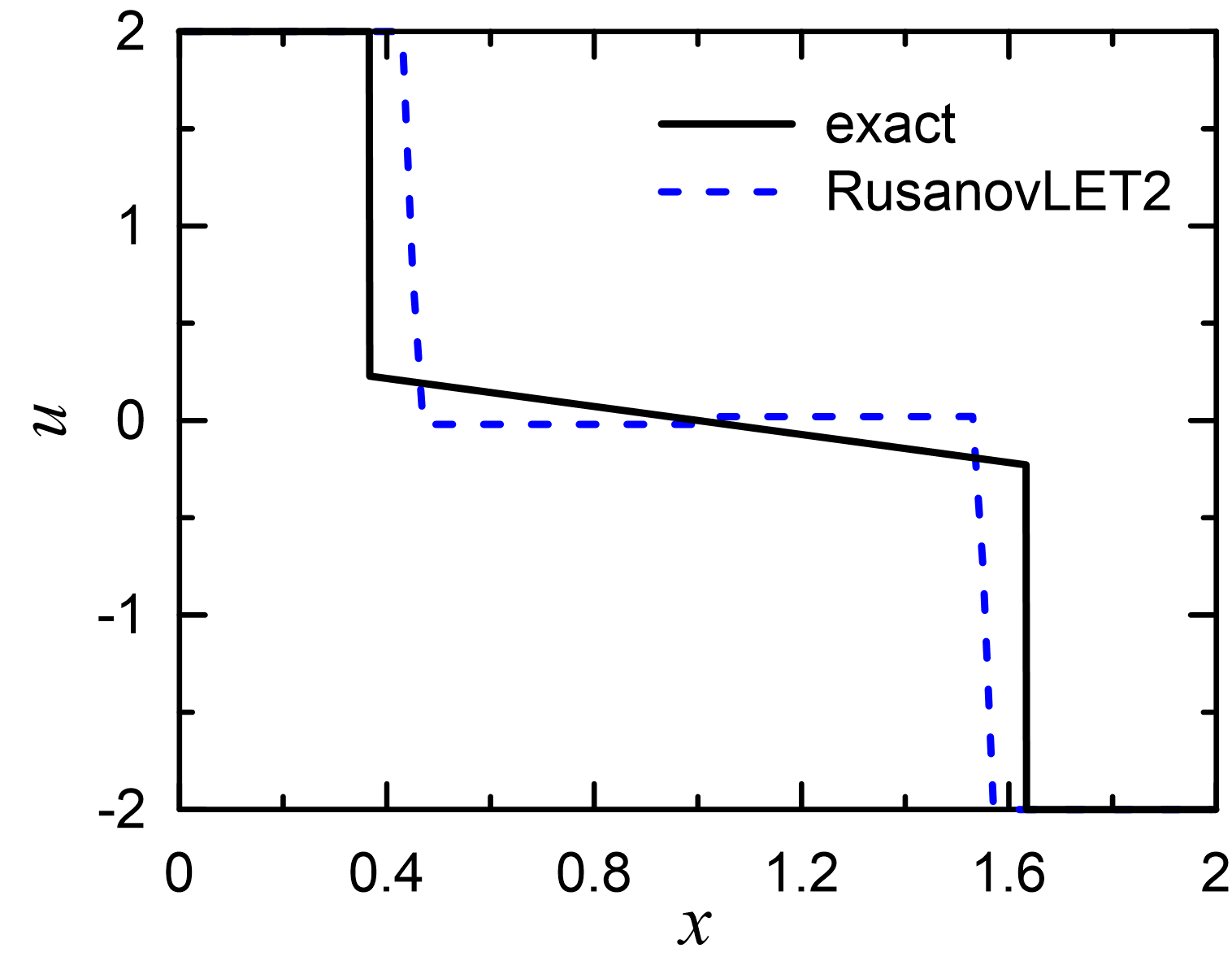}
\caption{Numerical solutions of the Riemann IVP for the scalar nonconvex conservation law  \eqref{eq:54}-\eqref{eq:56} obtained by applying the hybrid scheme \eqref{eq:513} with the second-order flux \eqref{eq:59}. Flux limiters on the left were calculated using linear programming without taking into account the discrete entropy inequality \eqref{eq:6}. Flux limiters on the right were calculated using linear programming taking into account the discrete cell entropy inequality \eqref{eq:514} with the Tadmor's numerical entropy flux \eqref{eq:57}}
\label{fig:3}       
\end{figure*}

According to Tadmor \cite[p.~464 Eq.~\eqref{eq:39}]{b12}, we approximate the entropy inequality \eqref{eq:3} by the following semi-discrete scheme
\begin{multline}
\label{eq:58}
   \frac{d}{{dt}}U({y_i}) + \frac{1}{{\Delta x}}\left[ {{H_{i + {1/2}}^T} - {H_{i - {1/2}}^T}} \right] =  \frac{1}{{2\Delta x}}\left( {{h_{i + {1 / 2}}}\,{\Delta _{i + {1 / 2}}}v\, - {\Delta _{i + {1 / 2}}}\psi } \right)  \\
   + \frac{1}{{2\Delta x}}\left( {{h_{i - {1 / 2}}}\,{\Delta _{i - {1 / 2}}}v\, - {\Delta _{i - {1 / 2}}}\psi } \right) \le 0  
\end{multline} 
where ${\Delta _{i + {1 / 2}}}v = {v_{i + 1}} - {v_i}$

In the numerical simulations, we will use fluxes of the second and fourth-order of spatial accuracy as high-order fluxes
\begin{eqnarray}
\label{eq:59}
 && h_{i + {1 / 2}}^{H2} = \frac{1}{2}({f_i} + {f_{i + 1}}) 										\\
 && h_{i + {1 / 2}}^{H4} = \frac{7}{{12}}({f_i} + {f_{i + 1}}) - \frac{1}{{12}}({f_{i + 2}} + {f_{i - 1}})								
\label{eq:510}
\end{eqnarray}

For the second-order flux $h_{i + {1 / 2}}^{H2}$ and the Rusanov numerical flux \eqref{eq:51}, the constraints \eqref{eq:32}-\eqref{eq:34} can be written in the form 
\begin{eqnarray}
\label{eq:511}
&& \frac{{\Delta x}}{{\Delta t}}\left( {\mathop {\min }\limits_{j \in {S_i}} {y_j} - {y_i}} \right) + h_{i + {1/2}}^{Rus} - h_{i - {1/2}}^{Rus} \le \nonumber \\
&& \le  - \frac{\alpha _{i + {1 / 2}}}{2}\mathop {\max }\limits_{s \in \left[ {{y_i},{y_{i + 1}}} \right]} \left| {f'(s)} \right|\;{\Delta _{i + {1/2}}}y + \frac{\alpha _{i - {1/2}}}{2}\mathop {\max }\limits_{s \in \left[ {{y_{i - 1}},{y_i}} \right]} \left| {f'(s)} \right|\;{\Delta _{i - {1 / 2}}}y \le \nonumber \\					
&& \le \frac{\Delta x}{\Delta t}\left( {\mathop {\max }\limits_{j \in {S_i}} {y_j} - {y_i}} \right) + h_{i + {1 / 2}}^{Rus} - h_{i - {1 / 2}}^{Rus}
\end{eqnarray} 

\begin{eqnarray}
\label{eq:512}
\begin{split}
  \frac{{\Delta x}}{{\Delta t}}\left[ {U({{\hat y}_i}) - U({y_i}) - {b_i}\left( {{{\hat y}_i} - {y_i}} \right)} \right]  + H_{i + {1/2}}^{Rus} - {b_i}h_{i + {1/2}}^{Rus} - H_{i - {1/2}}^{Rus} + {b_i}h_{i - {1/2}}^{Rus} \le \\				
  \le \frac{\alpha _{i + {1 / 2}}}{2}\mathop {\max }\limits_{s \in \left[ {{y_i},{y_{i + 1}}} \right]} \left| {f'(s)} \right|({b_i}{\Delta _{i + {1/2}}}y - {\Delta _{i + {1/2}}}U) - \\ 
  -{}\frac{\alpha _{i - {1/2}}}{2}\mathop {\max }\limits_{s \in \left[ {{y_{i - 1}},{y_i}} \right]} \left| {f'(s)} \right|({b_i}{\Delta _{i - {1/2}}}y - {\Delta _{i - {1/2}}}U)
\end{split}  			
\end{eqnarray} 

\begin{eqnarray}
\label{eq:513}
\begin{split}
 \frac{{\Delta x}}{{\Delta t}}\left( {{{\hat y}_i} - {y_i}} \right) + h_{i + {1 / 2}}^{Rus} + \frac{\alpha _{i + {1 / 2}}}{2}\mathop {\max }\limits_{s \in [{y_i},{y_{i + 1}}]} \left| {f'(s)} \right|{\Delta _{i + {1 / 2}}}y - \\
 - h_{i - {1/2}}^{Rus} - \frac{\alpha _{i - {1/2}}}{2}\mathop {\max }\limits_{s \in [{y_{i - 1}},{y_i}]} \left| {f'(s)} \right|{\Delta _{i - {1/2}}}y = 0 	
\end{split}  			
\end{eqnarray} 

Accordingly, for the Tadmor's numerical entropy flux, the cell entropy inequalities \eqref{eq:58}  can be written as
\begin{eqnarray}
\label{eq:514}
\begin{split}
 & {\alpha _{i + {1 / 2}}}\frac{{{\Delta _{i + {1 / 2}}}v}}{4}\;\mathop {\max }\limits_{s \in \left[ {{y_i},{y_{i + 1}}} \right]} \left| {f'(s)} \right|\;{\Delta _{i + {1 / 2}}} y \; + \\ 
 & + {\alpha _{i - {1 / 2}}}\frac{{{\Delta _{i - {1 / 2}}}v}}{4}\;\mathop {\max }\limits_{s \in \left[ {{y_{i - 1}},{y_i}} \right]} \left| {f'(s)} \right|\;{\Delta _{i - {1 / 2}}}y \le \\ 			 &  \le \frac{1}{2}\left( {{\Delta _{i + {1 / 2}}}\psi  - \,h_{i + {1/2}}^{Rus}\,{\Delta _{i + {1/2}}}v} \right) + \frac{1}{2}\left( {{\Delta _{i - {1/2}}}\psi  - h_{i - {1/2}}^{Rus}\,{\Delta _{i - {1/2}}}v\,} \right)  + \\
 & + \frac{{\Delta x}}{{\Delta t}}\left[ {{v_i}({{\hat y}_i} - {y_i}) - ({{\hat U}_i} - {U_i})} \right] 
\end{split}  			
\end{eqnarray} 

%
\begin{figure*}[!tb]
  \includegraphics[width=0.5\textwidth]{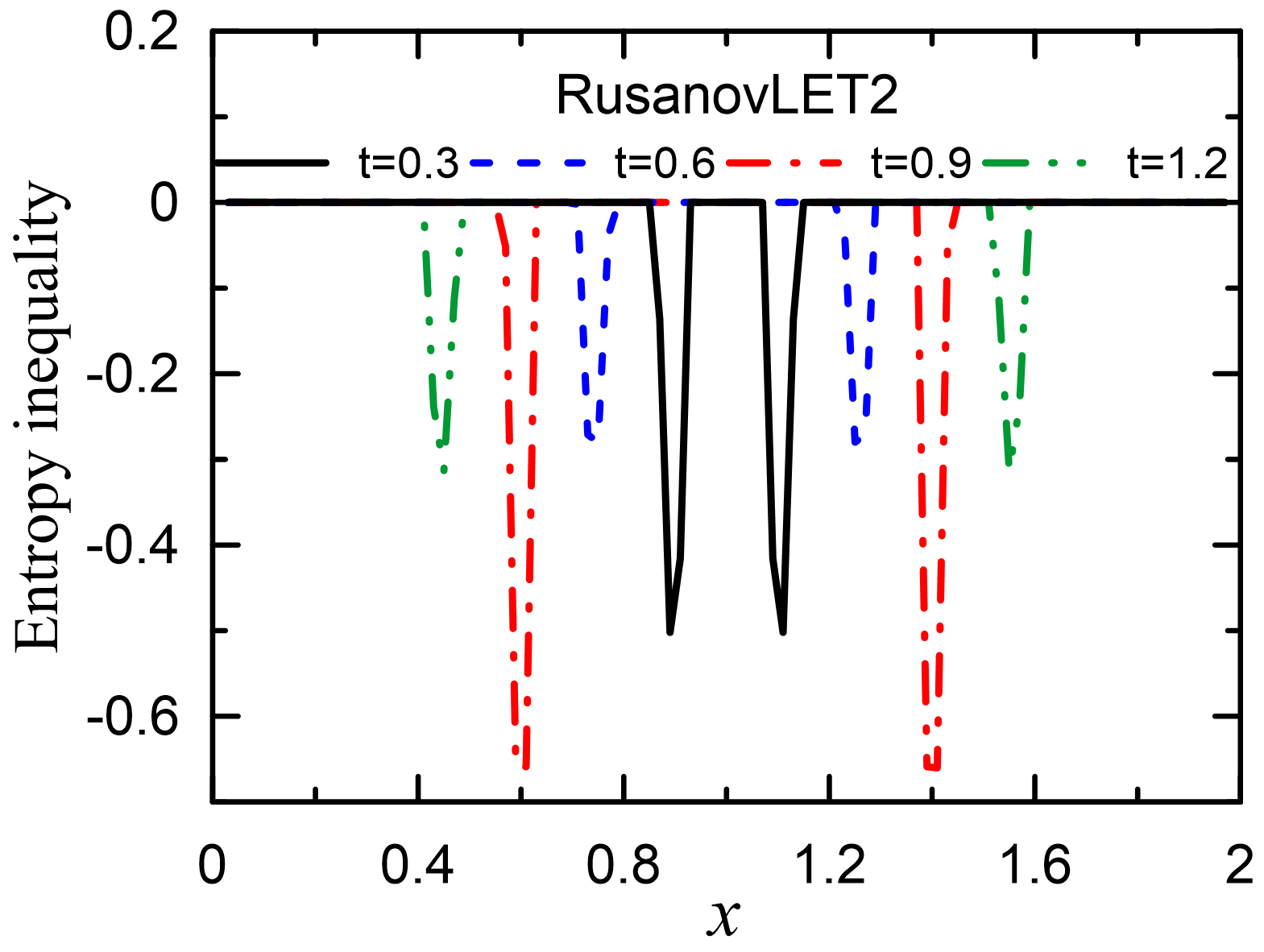}
  \includegraphics[width=0.5\textwidth]{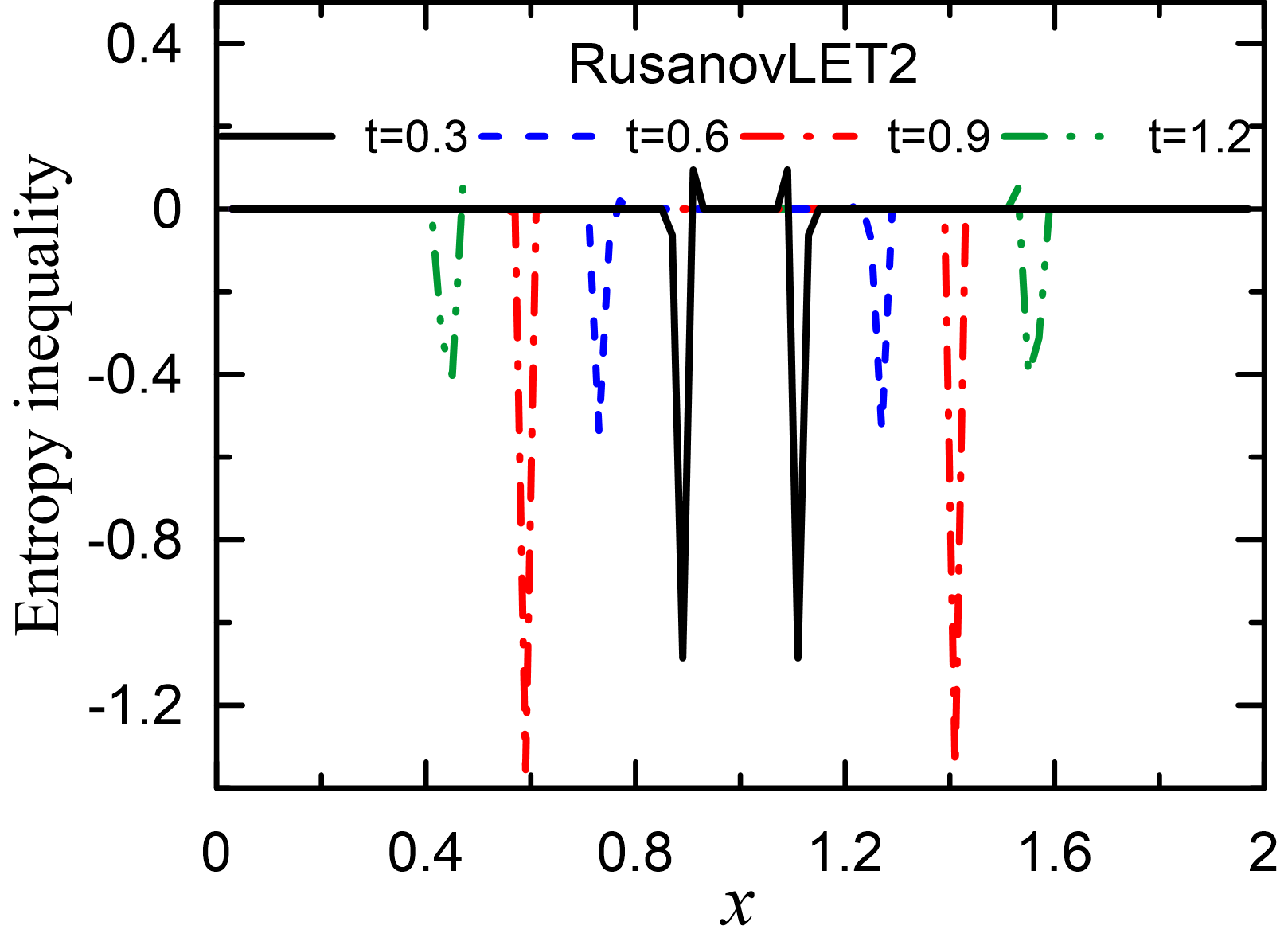}
\caption{Values of the discrete cell entropy inequality with the Tadmor's (left) and the proper (right) numerical entropy fluxes for numerical solution RusanovLET2 at different times. RusanovLET2 was obtained by applying the hybrid scheme with the second-order numerical flux \eqref{eq:59}, and flux limiters were calculated by using linear programming taking into account the discrete entropy inequalities \eqref{eq:514} with the Tadmor's numerical entropy flux}
\label{fig:4}       
\end{figure*}

We will mark a numerical solution with a scheme name and an ending that indicates how the flux limiters $\boldsymbol {\alpha}$ were calculated. The endings $LP$ and $LE$ mean that the flux limiters were calculated using linear programming without/with taking into account the cell entropy inequality \eqref{eq:6}, respectively. $AP$ and $AE$ mean that the flux limiters were calculated by using approximate relations. The letter $T$ indicates that the Tadmor's numerical entropy flux \eqref{eq:57} was applied to discretize the entropy inequality \eqref{eq:3}. The number at the end means  the order of spatial accuracy  of the numerical flux $h_{i + {1/2}}^H$. 

The computational domain is [0,2] with the space step of 0.02 and the temporal step of 0.002. We choose the square entropy function $U=0.5u^2$. Then the entropy flux $F$ and the entropy flux potential $\psi$ are defined by
\begin{eqnarray*}
 && F(u) = \left( \frac{u^2}{5}  - \frac{5}{6} \right) {u^3} \\ 
 && \psi (u) = u\,f - F = \frac{1}{{20}}{u^5} - \frac{5}{{12}}{u^3} + u
\end{eqnarray*} 
All numerical results are depicted at time $T=1.2$.

Numerical solutions of the Riemann IVP \eqref{eq:54}-\eqref{eq:56} obtained with the Godunov and Rusanov schemes are presented in Fig.~\ref{fig:1}. As shown in Fig.~\ref{fig:2} the numerical solution obtained with the Rusanov scheme satisfies the cell entropy inequality \eqref{eq:6} with the Tadmor's (left) and the proper (right) numerical entropy fluxes. 
%
\begin{figure*}[!t]
  \includegraphics[width=0.5\textwidth]{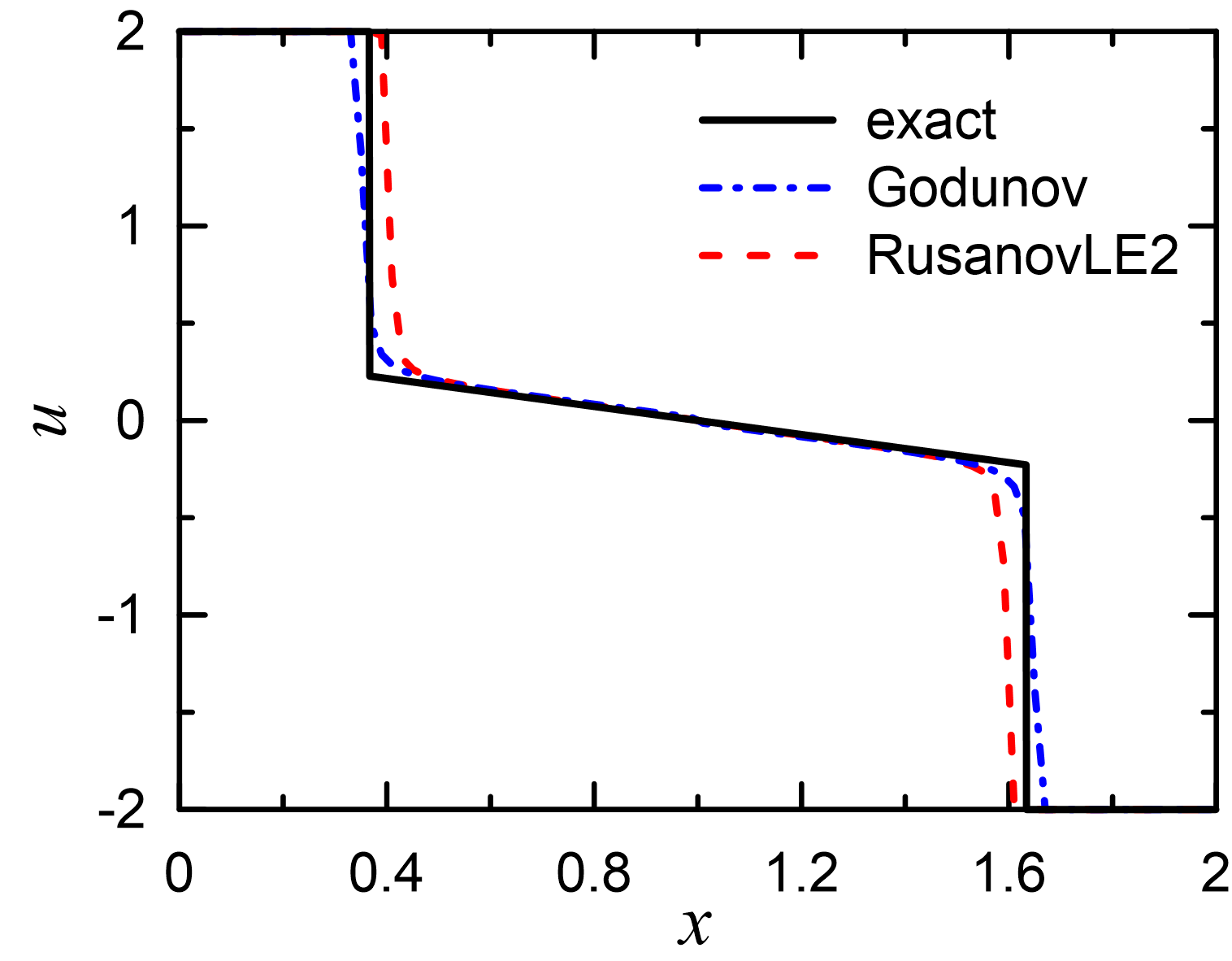}
  \includegraphics[width=0.5\textwidth]{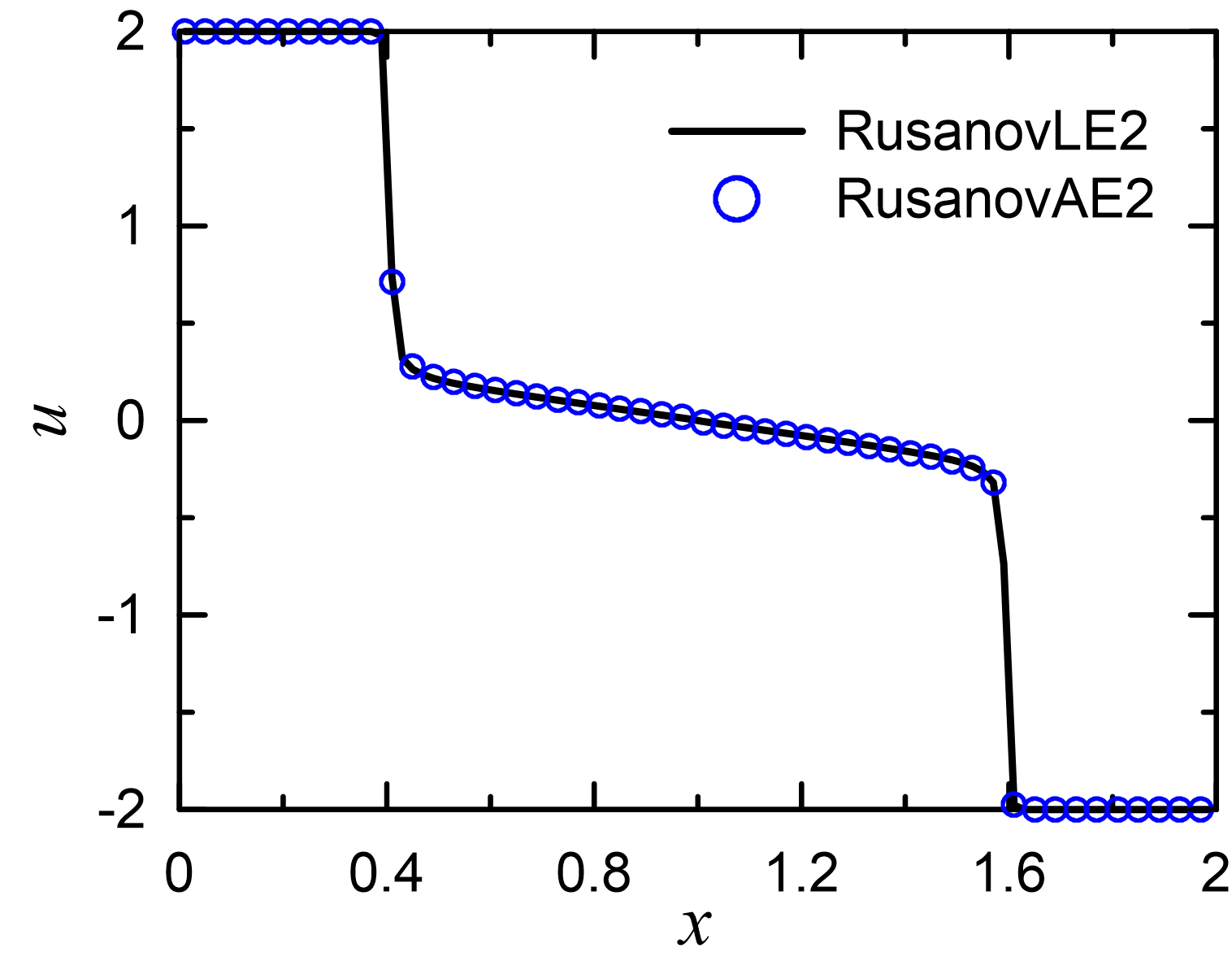}
\caption{Comparison of the numerical solutions Godunov, RusanovLE2 and RusanovAE2 for the Riemann IVP for the scalar nonconvex conservation law \eqref{eq:54}-\eqref{eq:56}}
\label{fig:5}       
\end{figure*}
%
%
\begin{figure*}[!t]
  \includegraphics[width=0.5\textwidth]{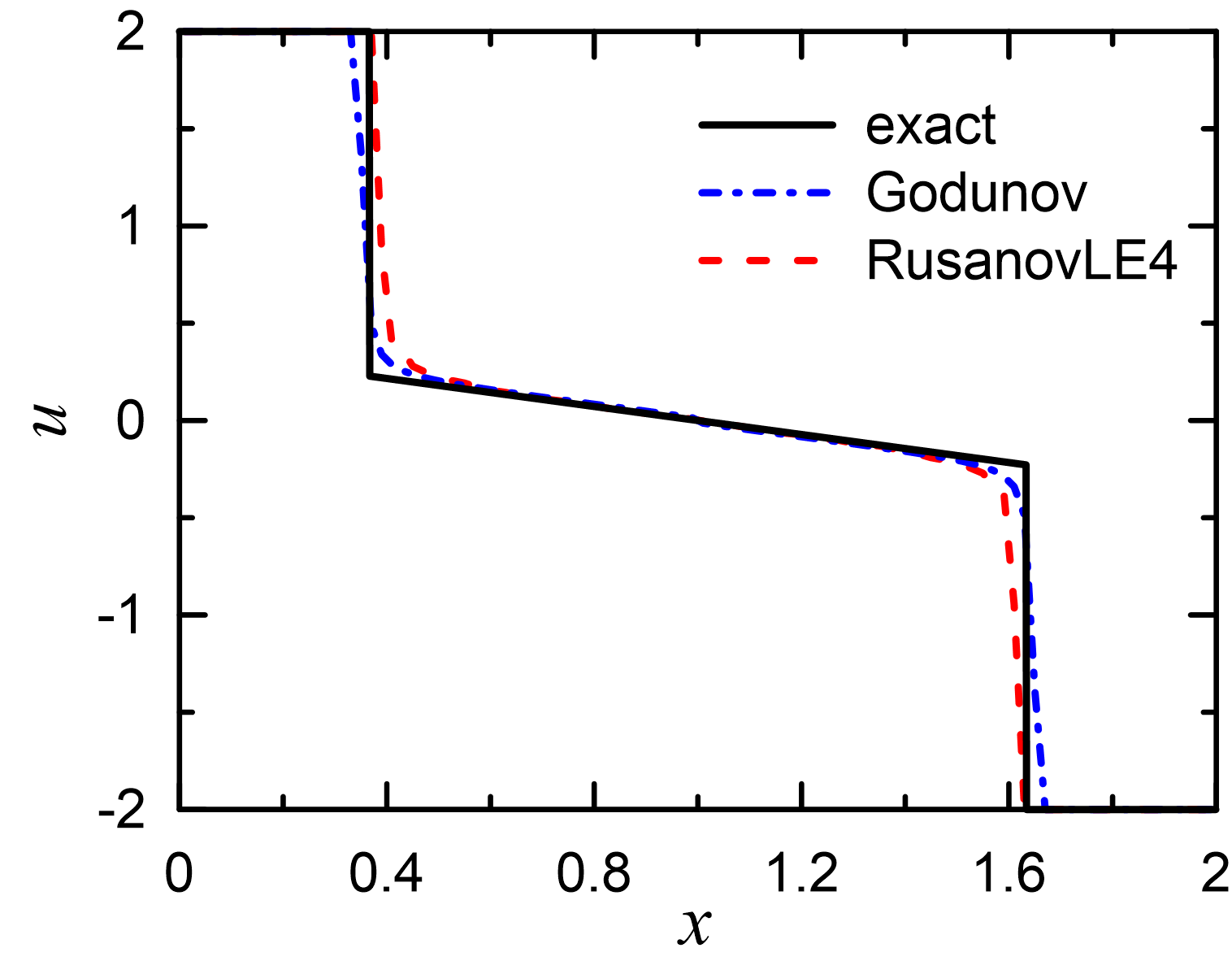}
  \includegraphics[width=0.5\textwidth]{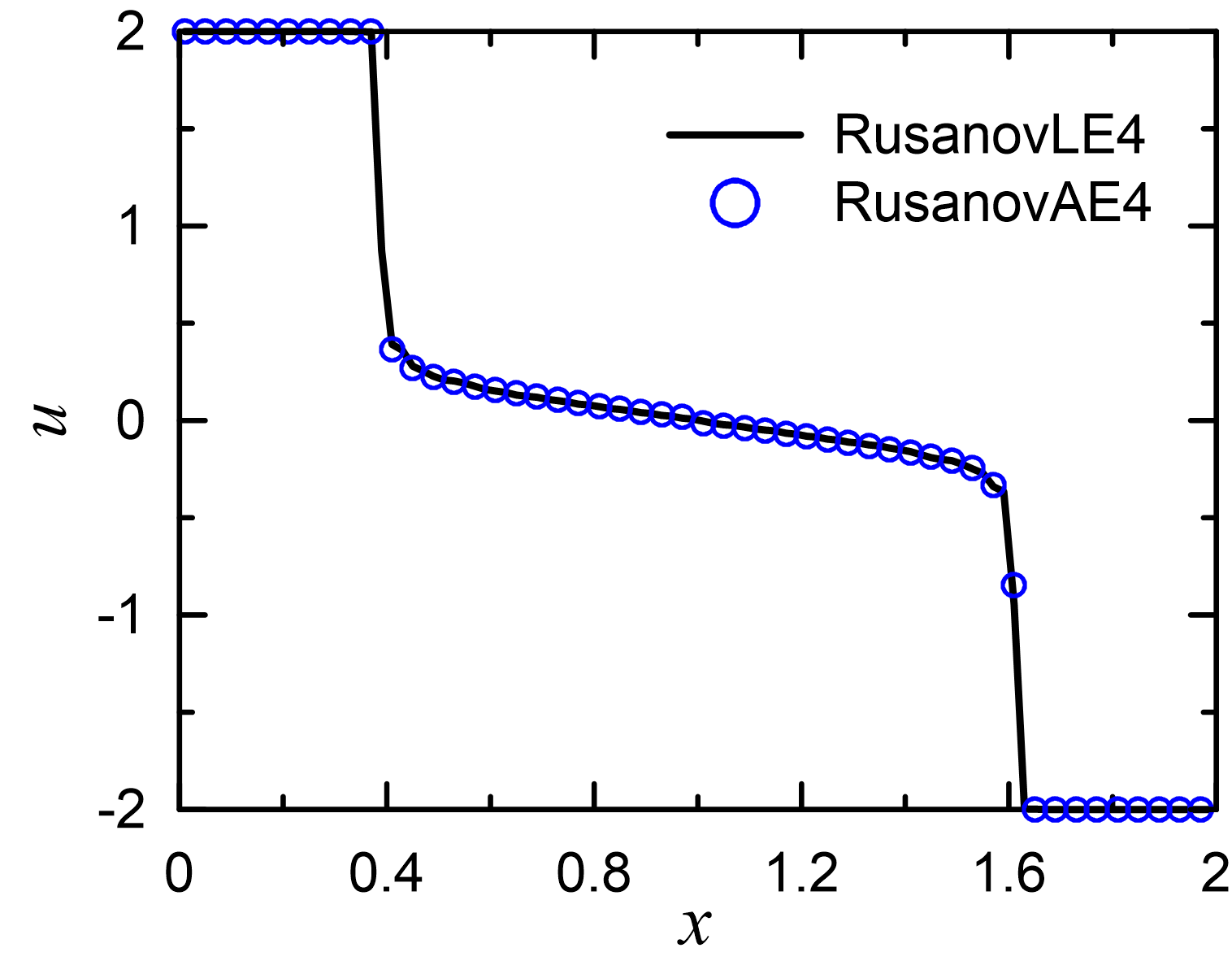}
\caption{Comparison of the numerical solutions Godunov, RusanovLE4 and RusanovAE4 for the Riemann IVP for the scalar nonconvex conservation law  \eqref{eq:54}-\eqref{eq:56}}
\label{fig:6}       
\end{figure*}

Simulation results of the Riemann IVP \eqref{eq:54}-\eqref{eq:56} obtained by the hybrid scheme \eqref{eq:513} with second-order flux \eqref{eq:59} are given in Fig.~\ref{fig:3}-Fig.~\ref{fig:5}. On the left Fig.~\ref{fig:3}  shows the numerical solution (RusanovLP2) for which the flux limiters were calculated by using linear programming without taking into account the cell entropy inequality \eqref{eq:6}. On the right Fig.~\ref{fig:3} presents the numerical solution (RusanovLET2) for which the flux limiters were calculated by linear programming taking into account the cell entropy inequality \eqref{eq:514} with the Tadmor's numerical entropy flux \eqref{eq:57}. One can see that the numerical solutions RusanovLP2 and RusanovLET2 are not physically relevant solutions of the Riemann IVP \eqref{eq:54}-\eqref{eq:56}. As shown in Fig.~\ref{fig:4}, RusanovLET2 satisfies the discrete cell entropy inequality with the Tadmor's numerical entropy flux (left) and violates it with the proper numerical entropy flux (right). Comparison of the numerical solutions Godunov, RusanovLE2 and Rusanov AE2 are presented in Fig.~\ref{fig:5}. We note a good agreement between the numerical solutions RusanovLE2 and RusanovAE2, for which the flux limiters were calculated by using exact and approximate solutions of the linear programming problem with the constraints \eqref{eq:511}-\eqref{eq:512}, respectively.
 
Comparison of the numerical solutions RusanovLE4 and RusanovAE4 with fourth-order accuracy in space for the Riemann IVP \eqref{eq:54}-\eqref{eq:56} is given in Fig.~\ref{fig:6}.

\subsection{Buckley-Leverett Equation}
\label{sec:52}

We consider the following Riemann problem for the Buckley-Leverett equation \cite{b35}
\begin{equation}
\label{eq:518a} 
  \frac{\partial u }{\partial t} + \frac{\partial }{\partial x}\left( {\frac{4u^2}{4u^2 + (1 - u )^2}} \right) = 0 
\end{equation}
subject to
\begin{equation}
\label{eq:519a} 
 u (x,0) = \begin{cases}
-3 \qquad & \rm{if} \;\; {\it x} < 0 \\
3  \qquad & \rm{if} \;\; {\it x} \ge 0
\end{cases}  
\end{equation}

The exact entropy solution consists of two shock waves and a rarefaction wave that is close to 0. In~\cite{b35} Chen and Shu developed the entropy stable high-order discontinuous Galerkin method based on the cell entropy inequalities with the Tadmor's numerical entropy fluxes. They applied their scheme to obtain the numerical solution of the Riemann problem \eqref{eq:518a}-\eqref{eq:519a} for the square entropy function and obtained a physically noncorrect solution. Therefore, our purpose is to repeat this test with the same initial data.

%
\begin{figure*}[!tb]
  \includegraphics[width=0.5\textwidth]{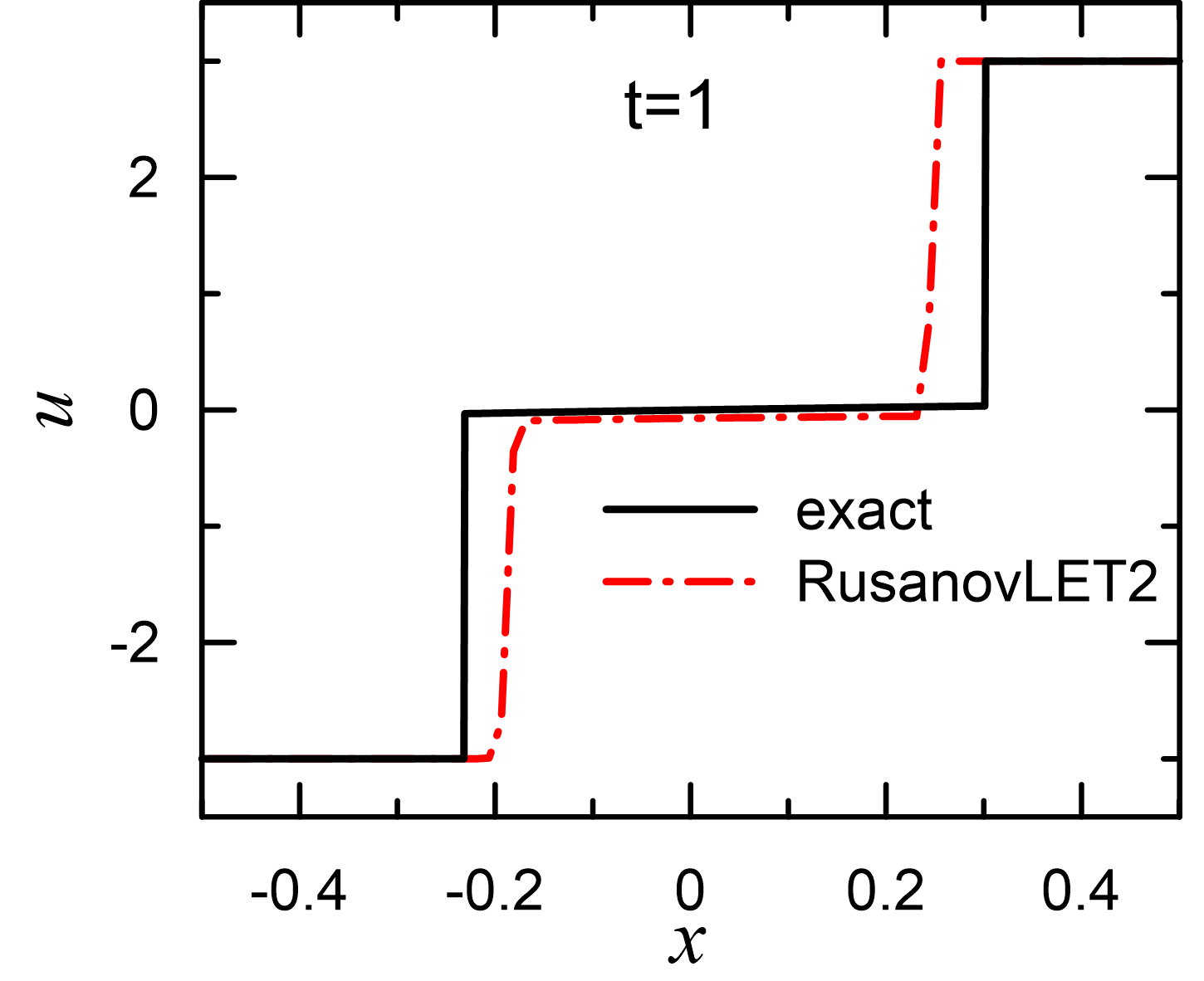}
  \includegraphics[width=0.5\textwidth]{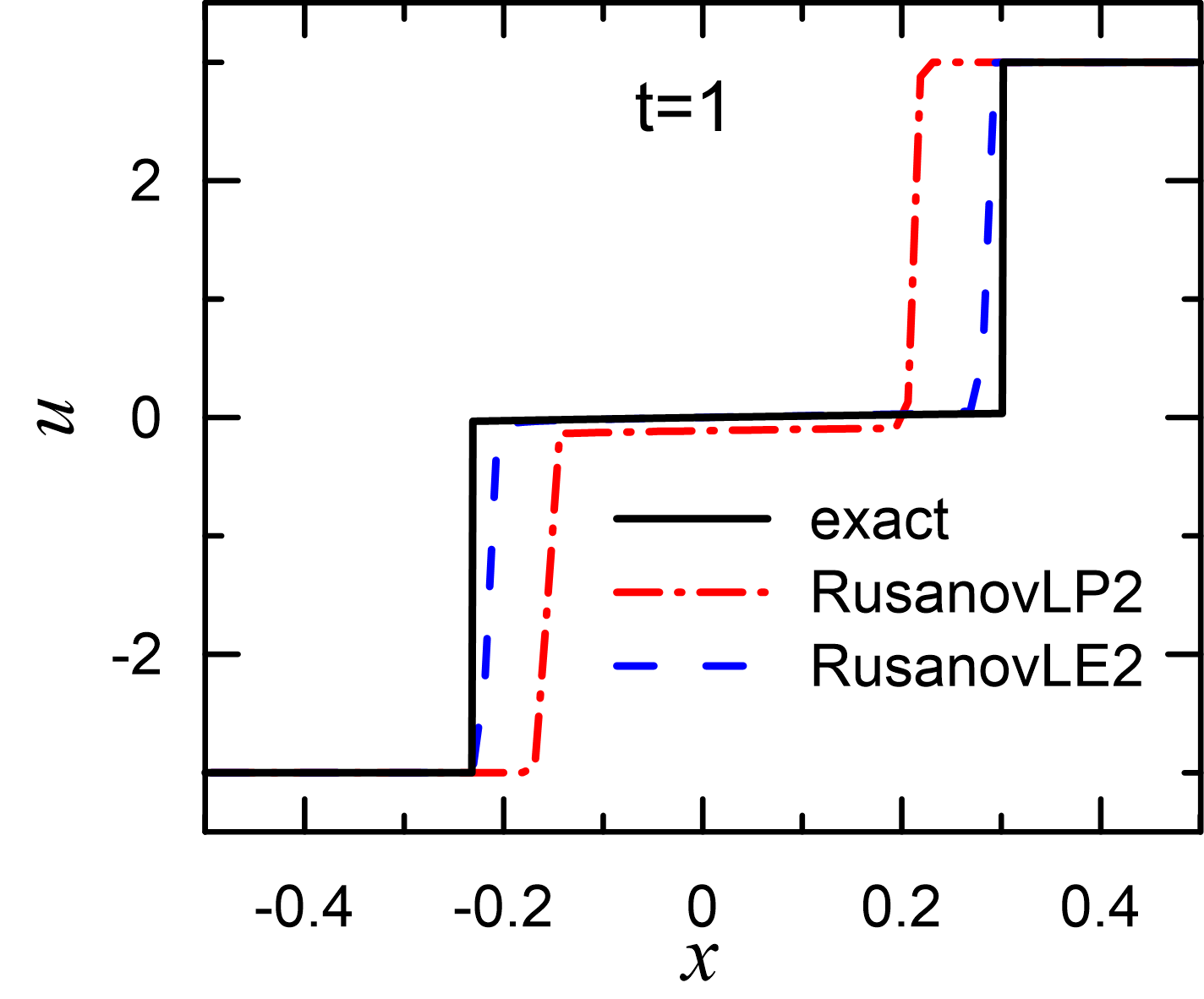}
\caption{Numerical solutions of the Riemann problem of the Buckley–Leverett equation at $t=1$ by using the hybrid scheme \eqref{eq:513}. Flux limiters are computed by using linear programming with and without taking into account the discrete entropy inequality. On the left the numerical solution that was obtained by using the discrete cell entropy inequality \eqref{eq:514} with the Tadmor's numerical entropy flux \eqref{eq:57} and on the right the numerical solution obtained with the cell entropy inequality with the proper numerical entropy flux}
\label{fig:7}       
\end{figure*}
%

%
\begin{figure*}[!tb]
  \includegraphics[width=0.5\textwidth]{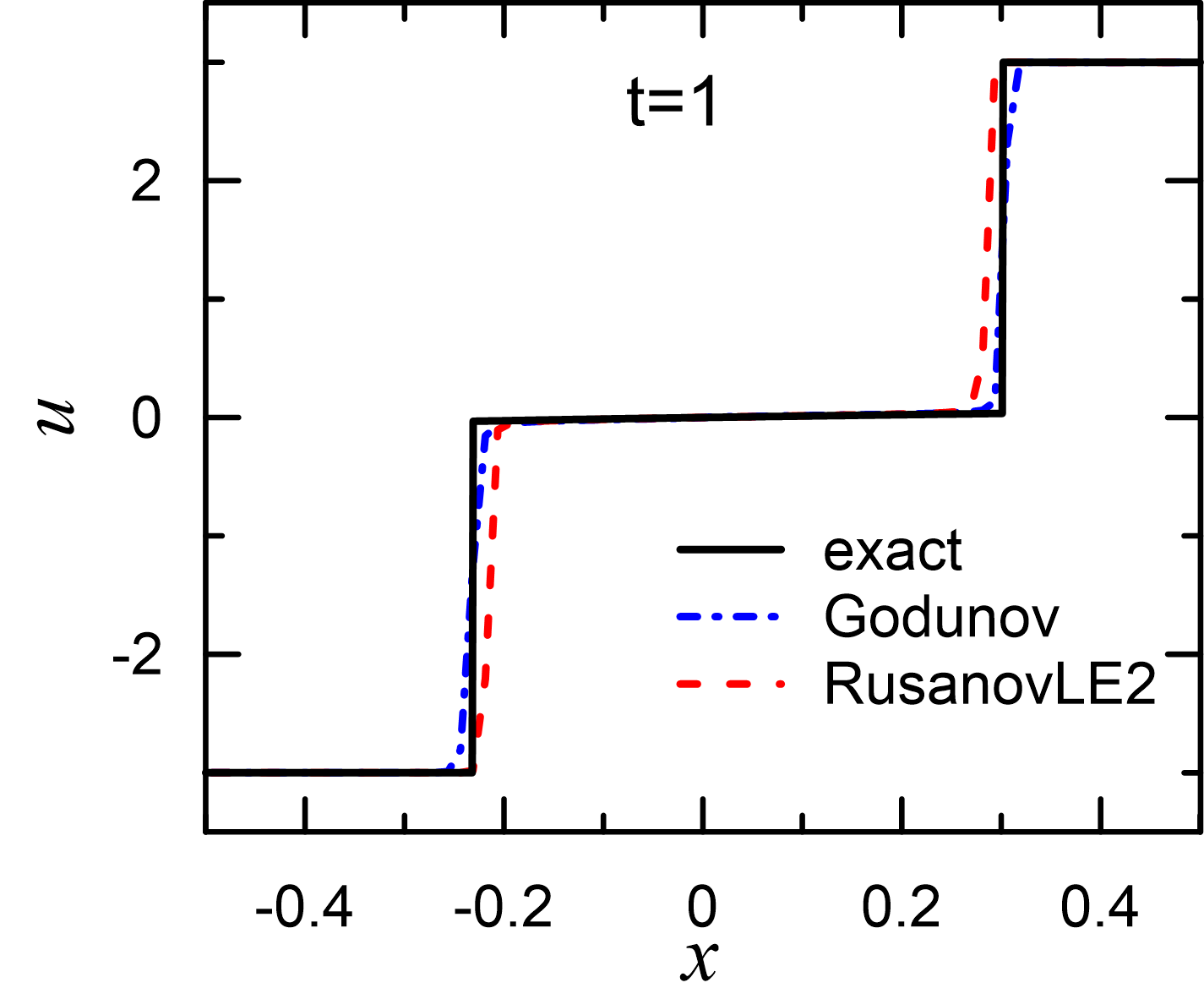}
  \includegraphics[width=0.5\textwidth]{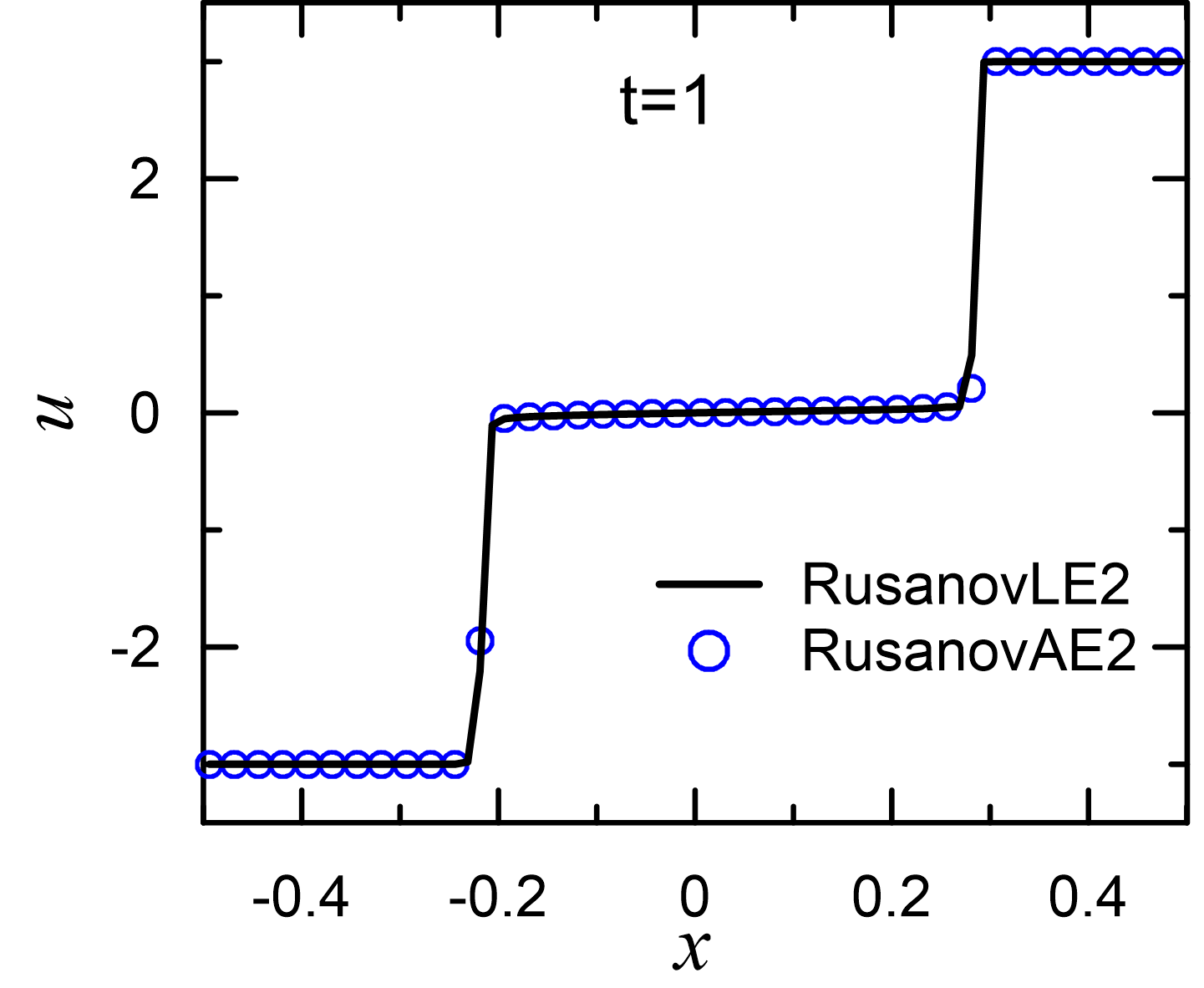}
\caption{Comparison between the numerical solutions of the Riemann problem of the Buckley-Leverett equation obtained with the Godunov and the hybrid scheme \eqref{eq:513}. Flux limiters are computed by using the exact (RusanovLE) and approximate (RusanovAE) solutions of the linear programming problem by taking into account the discrete entropy inequality \eqref{eq:512} }
\label{fig:8}       
\end{figure*}

Similarly to \cite{b35}, the computational domain is [-0.5,0.5] and the end time $t=1$. The entropy flux corresponding to the square entropy function   is defined as

  \[ F =  - \frac{4}{{25}}\left[ {\frac{{u + 2}}{{5{u^2} - 2u + 1}} + \ln \left( {5{u^2} - 2u + 1} \right) - \frac{3}{2}{{\tan }^{ - 1}}\left( {\frac{{5u - 1}}{2}} \right)} \right]   \]

Our numerical solutions obtaining with the hybrid scheme \eqref{eq:513} on 80 cells are given in Fig.~\ref{fig:7}. The flux limiters are calculated by using linear programming with and without taking into account discrete entropy inequalities \eqref{eq:512} and \eqref{eq:514}. Obviously, the numerical solution RusanovLP2 that is obtained without taking into account the discrete entropy inequality \eqref{eq:512} and corresponding to the classical FCT solution is not physically correct. Also, the numerical solution RusanovLET2 obtained by using the discrete entropy inequality with the Tadmor's numerical entropy flux is not physically relevant.

Comparison between the numerical solutions RusanovLE2 and RusanovAE2 is shown in Fig.~\ref{fig:8}. The flux limiters for RusanovLE2 and RusanovAE2 are computed by using the exact and approximate solutions of linear programming problems with the constraints \eqref{eq:511}-\eqref{eq:512}. We note a good agreement of these numerical solutions.

\section{Conclusions}
\label{sec:5}

We presented the flux correction design for a hybrid scheme to obtain a physically relevant solution of scalar hyperbolic conservation laws. The hybrid scheme is a linear combination of a monotone scheme and a high-order scheme. The flux limiters for the hybrid scheme are a solution of a corresponding optimization problem in which constraints valid for the monotone scheme are applied to the hybrid scheme. It is shown that the physically relevant solution for scalar hyperbolic conservation laws satisfies the discrete entropy inequality with the proper numerical entropy flux. A numerical solution that satisfies the discrete entropy inequality with the Tadmor's numerical entropy flux maybe not a physically correct. An approximate solution of the optimization problem can be chosen to compute flux limiters. Numerical examples show a good agreement between numerical solutions of scalar hyperbolic conservation laws for which the flux limiters are calculated from exact and approximate solutions of the optimization problem.


%
%



\end{document}